\documentclass[10pt,english,leqno]{article}
\usepackage[T1]{fontenc}
\usepackage{times}
\usepackage[latin9]{inputenc}
\usepackage{amsthm}
\usepackage{amsmath}
\usepackage{amssymb}
\usepackage{enumerate}
\usepackage{hyperref}
\usepackage{tkz-graph}
\usepackage{todonotes}

\makeatletter

\newtheorem{thm}{Theorem}
\newtheorem{lemma}[thm]{Lemma}
\newtheorem{prop}[thm]{Proposition}
\newtheorem{cor}[thm]{Corollary}
\theoremstyle{definition}

\DeclareMathOperator{\rank}{rank}

\DeclareMathOperator{\Tr}{Tr}

\newsavebox{\fmbox}

\newcommand{\oR}{\ensuremath{\mathbb{R}}}
\newcommand{\oN}{\ensuremath{\mathbb{N}}}

\newlength{\dyindent}
{\end{list}}

\newenvironment{dy*}{\refstepcounter{equation}\begin{list}{}%
{\setlength{\leftmargin}{\dyindent}\setlength{\labelwidth}{\dyindent}%
\addtolength{\labelwidth}{-\labelsep}}%
\item}%
{\end{list}}

\newcommand{\mr}{\mbox{mr}}

\newcommand{\nullity}{\mbox{nullity}}
\newcommand{\Is}{\mathcal{I}^s}
\newcommand{\I}{\mathcal{I}}
\newcommand{\pin}{\mbox{pin}}

\title{The inertia set of a signed graph}

\author{Marina Arav, Frank J. Hall, Zhongshan Li, Hein van der Holst\footnote{Corresponding author, E-mail: hvanderholst@gsu.edu} \\
Department of Mathematics and Statistics \\
Georgia State University \\
Atlanta, GA 30303, USA
}
\date{}

\begin{document}

\maketitle

\begin{abstract}
A signed graph is a pair $(G,\Sigma)$, where $G=(V,E)$ is a graph (in which parallel edges are permitted, but loops are not) with $V=\{1,\ldots,n\}$ and $\Sigma\subseteq E$. By $S(G,\Sigma)$ we denote the set of all symmetric $V\times V$ matrices $A=[a_{i,j}]$ with $a_{i,j}<0$ if $i$ and $j$ are connected by only even edges, $a_{i,j}>0$ if $i$ and $j$ are connected by only odd edges, $a_{i,j}\in \mathbb{R}$ if $i$ and $j$ are connected by both even and odd edges, $a_{i,j}=0$ if $i\not=j$ and $i$ and $j$ are non-adjacent, and $a_{i,i} \in \mathbb{R}$ for all vertices $i$. The stable inertia set of a signed graph $(G,\Sigma)$ is the set of all pairs $(p,q)$ for which there exists a matrix $A\in S(G,\Sigma)$ with $p$ positive and $q$ negative eigenvalues which has the Strong Arnold Property. In this paper, we study the stable inertia set of (signed) graphs.
\end{abstract}

\noindent keywords: graph, signed graph, inertia, symmetric, minor\newline
MSC: 05C05, 05C22, 05C50, 05C83, 15A03

\newpage

\section{Introduction}

A \emph{signed graph} is a pair $(G,\Sigma)$, where $G=(V,E)$ is a graph (in which parallel edges are permitted, but loops are not)  with $V=\{1,\ldots,n\}$ and $\Sigma\subseteq E$. (We refer to \cite{Diestel} for the notions and concepts in Graph Theory.) The edges in $\Sigma$ are called \emph{odd} and the other edges of $E$ \emph{even}. If $V=\{1,2,\ldots,n\}$, we denote by $S(G,\Sigma)$ the set of all real symmetric $n\times n$ matrices $A=[a_{i,j}]$ with 
\begin{itemize}
\item $a_{i,j} < 0$ if $i$ and $j$ are connected by only even edges, \item $a_{i,j}>0$ if $i$ and $j$ are connected by only odd edges, 
\item $a_{i,j}\in \mathbb{R}$ if $i$ and $j$ are connected by both even and odd edges, 
\item $a_{i,j}=0$ if $i\not=j$ and $i$ and $j$ are non-adjacent, and 
\item $a_{i,i} \in \mathbb{R}$ for all vertices $i$. 
\end{itemize}
If $A$ is a symmetric matrix, then by $\pin(A)$ we denote the pair $(p,q)$, where $p$ and $q$ are the number of positive and negative eigenvalues of $A$, respectively. We define the \emph{inertia set} of a signed graph $(G,\Sigma)$ as the set $\{\pin(A)~|~A\in S(G,\Sigma)\}$ and denote it by $\I(G,\Sigma)$. The analogous version for graphs was introduced by Barrett, Hall, and Loewy in \cite{MR2547901}. For a graph $G$, denote by $S(G)$  the set of all real symmetric $n\times n$ matrices $A=[a_{i,j}]$ with $a_{i,j}\not=0$ if $i$ and $j$ are connected by a single edge, $a_{i,j}\in \mathbb{R}$ if $i$ and $j$ are connected by multiple edges, $a_{i,j}=0$ if $i\not=j$ and $i$ and $j$ are non-adjacent, and $a_{i,i} \in \mathbb{R}$ for all vertices $i$.  The inertia set of a graph $G$ is the set $\{\pin(A)~|~A\in S(G)\}$ and is denoted by $\I(G)$. The inertia set of a signed graph is a refinement of the inertia set of a graph: If $G=(V,E)$ is a graph, then $\I(G) = \cup_{\Sigma\subseteq E} \I(G,\Sigma)$. 

We also make the following definitions for signed graphs. The \emph{minimum rank} of a signed graph $(G,\Sigma)$, denoted $\mr(G,\Sigma)$, is the minimum of the ranks of the matrices in $S(G,\Sigma)$. The \emph{minimum semidefinite rank} of a signed graph $(G,\Sigma)$, denoted $\mr_+(G,\Sigma)$, is the minimum of the ranks of the positive semidefinite matrices in $S(G,\Sigma)$. The \emph{maximum nullity} of a signed graph $(G,\Sigma)$, denoted $M(G,\Sigma)$, is the maximum of the nullities of the matrices in $S(G,\Sigma)$, and the \emph{maximum semidefinite nullity}, denoted $M_+(G,\Sigma)$, is the maximum of the nullities of the positive semidefinite matrices in $S(G,\Sigma)$. It is clear that if $(G,\Sigma)$ has $n$ vertices, then $M(G,\Sigma) + \mr(G,\Sigma) = n$ and $M_+(G,\Sigma) + \mr_+(G,\Sigma) = n$. 
The minimum rank, minimum semidefinite rank, the maximum nullity, and the maximum semidefinite nullity of a graph $G$ are defined in the same way as the corresponding parameters for signed graphs, except that one replaces $S(G,\Sigma)$ by $S(G)$. Clearly, for a graph $G=(V,E)$, $\mr(G) = \min_{\Sigma\subseteq E} \mr(G,\Sigma)$ and $\mr_+(G) = \min_{\Sigma\subseteq E}\mr_+(G,\Sigma)$.
The inertia set generalizes the minimum rank and minimum semidefinite rank. 
The minimum rank (minimum semidefinite rank) is equal to the smallest integer $k\geq 0$ such that there exists a pair $(p,q)\in \I(G,\Sigma)$ with $p+q=k$ (such that there exists a pair $(p,0)\in \I(G,\Sigma)$ with $p=k$). Part of the results in this paper are characterizations of the classes of signed graphs $(G,\Sigma)$ with $M_+(G,\Sigma)\leq 1$ and with $M(G,\Sigma)\leq 1$. 
Currently, the study of minimum rank and minimum semidefinite rank of a graph is an area of active research. For a survey on the minimum and minimum semidefinite ranks of graphs, we refer to \cite{FalHog2007}. 

For an integer $n\geq 0$, let $\mathbb{N}^2_{[0,n]}$ denote the set of all pair $(p,q) \in \mathbb{N}^2$ with $p+q\leq n$. (Here we include $0$ in the set of all natural numbers $\mathbb{N}$.) To obtain the inertia set of a signed graph $(G,\Sigma)$ with $n$ vertices, we need to check for each $(p,q)\in \mathbb{N}^2_{[0,n]}$ whether or not there exists a matrix $A\in S(G,\Sigma)$ with $p$ positive and $q$ negative eigenvalues. Some pairs $(p,q)$ for which there exists such a matrix $A\in S(G,\Sigma)$ can be found using the stable inertia set (defined below). In general, using the stable inertia set we do not find all pairs of the inertia set of a signed graph, but we will see that in many cases the stable inertia set is already sufficient to determine the inertia set of a signed graph, and if not, we can at least partially determine the inertia set of a signed graph.  Before describing what we mean by the stable inertia set of a signed graph, we describe graph parameters that are special cases of the stable inertia set.
 
Colin de Verdi\`ere introduced in \cite{CdeV1} the interesting graph parameter $\mu$. In order to describe this parameter we need the notion of \emph{Strong Arnold Property} (\emph{SAP} for short). A matrix $A=[a_{i,j}]\in S(G)$ has the SAP if $X=0$ is the only symmetric matrix $X=[x_{i,j}]$ such that $x_{i,j} = 0$ if $i$ and $j$ are adjacent vertices or $i=j$, and $A X = 0$. For a graph $G$, $\mu(G)$ is defined as the largest nullity of any matrix $A\in S(G,\emptyset)$ that has exactly one negative eigenvalue and has the SAP\@. This parameter has a very nice property: if $H$ is a minor of $G$, then $\mu(H)\leq \mu(G)$. Recall that a minor of a graph $G$ is any graph that can be obtained from $G$ by a series of contractions of edges of subgraph of $G$. Further, the parameter $\mu$ characterizes outerplanar and planar graphs as those graphs $G$ such that $\mu(G)\leq 2$ and $\mu(G)\leq 3$, respectively, see \cite{CdeV1}. Lov\'asz and Schrijver \cite{LovaszSchrijver98a} showed that graphs $G$ that have a flat embedding are exactly those for which $\mu(G)\leq 4$. For a survey on $\mu$, see \cite{HLovaszS2}.

In \cite{CdeV3}, Colin de Verdi\`ere introduced the graph parameter $\nu$. For a simple graph $G$, $\nu(G)$ is defined to be the largest nullity of any positive semidefinite matrix $A\in S(G)$ having the SAP\@. This parameter also has the very nice property that if $H$ is a minor of $G$, then $\nu(H)\leq \nu(G)$. (As $\nu(G)$ is a lower bound of $M_+(G)$, it can be used to obtain lower bounds for $M_+(G)$ using minors of $G$.) For the complete graph $K_n$ with $n>1$, $\nu(K_n)=n-1$. So if a simple graph $G$ contains a cycle, that is, if $G$ contains a minor isomorphic to $K_3$, then $\nu(K_3)=2\leq M_+(G)$. Hence a simple graph $G$ with $M_+(G)\leq 1$ does not contain any cycle, that is, $G$ is a forest. Colin de Verdi\`ere showed in \cite{CdeV3} that for simple graphs, $\nu(G)\leq 1$ if and only if $G$ is a forest. 
The simple graphs $G$ with $\nu(G)\leq 2$ have been characterized by Kotlov \cite{Kotlov2000a}. The parameter $\nu$ can be extended to graphs in which we allow parallel edges, but no loops; see \cite{Holst96a}. In \cite{Holst97a}, van der Holst gave a characterization of graphs $G$ with $\nu(G)\leq 2$, where parallel edges are permitted.

Barioli, Fallat, and Hogben introduced in \cite{BarFalHog2005a} a similar graph parameter $\xi$ for simple graphs $G$. This parameter is defined as the largest nullity of any matrix $A\in S(G)$ having the SAP\@. It has the same property that if $H$ is a minor of $G$, then $\xi(H)\leq \xi(G)$. For complete graphs $K_n$ with $n>1$, $\xi(K_n)=n-1$. Furthermore, $\xi(K_{1,3}) = 2$. So if a simple graph $G$ contains a cycle or a vertex with degree $>2$, then $2\leq \xi(G)\leq M(G)$. Only a disjoint union of paths satisfy the conditions of having no cycles and no vertices with degree $>2$. Thus, if $\xi(G)\leq 1$, then $G$ is a disjoint union of paths. Barioli, Fallat, and Hogben showed in \cite{BarFalHog2005a} that $\xi(G)\leq 1$ if and only if $G$ is a disjoint union of paths. The graphs $G$ with $\xi(G)\leq 2$ also have been characterized, see  \cite{HvdH2007a}.

We define the \emph{stable inertia set} of a signed graph $(G,\Sigma)$ as the set 
\begin{equation*}
\{\pin(A)~|~A\in S(G,\Sigma)\text{ and $A$ has the SAP}\},
\end{equation*}
and denote it by $\Is(G,\Sigma)$. It is clear that $\Is(G,\Sigma)\subseteq \I(G,\Sigma)$ for any signed graph $(G,\Sigma)$. Analogously, we define the stable inertia set of a graph $G=(V,E)$ as the set
\begin{equation*}
\{\pin(A)~|~A\in S(G)\text{ and $A$ has the SAP}\},
\end{equation*}
and denote it by $\Is(G)$. It is easy to see that $\Is(G) = \cup_{\Sigma\subseteq E} \Is(G,\Sigma)$ and that $\Is(G)\subseteq \I(G)$. 

If $v$ is a vertex of $(G,\Sigma)$, then $\delta(v)$ denotes the set of all edges incident with $v$. We call the operation $\Sigma\to \Sigma\Delta\delta(v)$ \emph{resigning around $v$}, where $\Delta$ is the symmetric difference. If $U\subseteq V$, then $\delta(U)$ denotes the set of all edges that have one end in $U$ and one end not in $U$. We say that $(G,\Sigma)$ and $(G,\Sigma\Delta\delta(U))$ are \emph{sign-equivalent} and call the operation $\Sigma\to \Sigma\Delta\delta(U)$ \emph{resigning on $U$}. Resigning on $U$ amounts to performing a diagonal similarity on the matrices in $S(G,\Sigma)$, and hence it does not affect the inertia set. We call a cycle $C$ of a signed graph $(G,\Sigma)$ \emph{odd} if $C$ has an odd number of odd edges, otherwise  we call $C$ \emph{even}. We call a signed graph \emph{bipartite} if it has no odd cycles. Zaslavsky showed in \cite{MR676405} that two signed graphs are sign-equivalent if and only if they have the same set of odd cycles.
Thus, signed graphs that have the same set of odd cycles have the same (stable) inertia set (and the same minimum rank and the same minimum semidefinite rank).

In Section~\ref{sec:minorsstableinertiaset}, the main section of the paper, we will show that the stable inertia set of a signed graph behaves well under taking subgraphs and contracting edges. (Contracting an edge $e$ with ends $u$ and $v$ means deleting $e$ and identifying $u$ and $v$, and since we do not allow loops, also deleting any loops that appear.) From this result it follows that the same holds for the stable inertia set of a graph. The stable inertia set of a signed graph includes in a simple way the graph parameters $\mu, \nu,$ and $\xi$. For a graph $G=(V,E)$ with $n$ vertices, $\mu(G)$ is the largest integer $k\geq 0$ such that $(n-k-1,1)\in \Is(G,\emptyset)$, $\nu(G)$ is the largest integer $k\geq 0$ such that $(n-k,0)\in \Is(G)$, and $\xi(G)$ is the largest integer $k\geq 0$ for which there exists an integer $p\geq 0$ such that $(p,n-k-p)\in \Is(G)$.

For a signed graph $(G,\Sigma)$, we define $\nu(G,\Sigma)$ as the largest nullity of any positive semidefinite matrix $A\in S(G,\Sigma)$ having the SAP, and $\xi(G,\Sigma)$ as the largest nullity of any matrix $A\in S(G,\Sigma)$ having the SAP\@. If $(G,\Sigma)$ has $n$ vertices, then $\nu(G,\Sigma)$ is the largest integer $k\geq 0$ such that $(n-k,0)\in \Is(G,\Sigma)$ and $\xi(G,\Sigma)$ is the largest integer $k\geq 0$ for which there exists an integer $p\geq 0$ such that $(p,n-k-p)\in \Is(G,\Sigma)$. A \emph{minor} of a signed graph $(G,\Sigma)$ is any signed graph that can be obtain from $(G,\Sigma)$ by deleting edges and vertices, contracting even edges, and resiging around vertices. From the result on the behavior of the stable inertia set of a signed graph under taking subgraphs and contracting edges, it follows that if $(H,\Omega)$ is a minor of $(G,\Sigma)$, then $\nu(H,\Omega)\leq \nu(G,\Sigma)$, and that if $H'$ is a minor of $G$, then $\xi(H',\Sigma\cap E(H'))\leq \xi(G,\Sigma)$.

Let us now introduce some notation. Let $A=[a_{i,j}]$ be an $n\times n$ matrix. If $R,S\subseteq \{1,\ldots,n\}$, we denote by $A[R,S]$ the submatrix of $A$ that lies in the rows of $A$ indexed by $R$ and the column indexed by $S$. If $R=S$, we will write $A[S]$ for $A[R,S]$; $A[\{v\}, S]$ can be denoted by $A[v,S]$. If $S\subseteq \{1,\ldots,n\}$, we let $\overline{S} = \{1,\ldots,n\}\setminus S$. If $A[S]$ is invertible, the Schur complement of $A[S]$ in $A$ is the matrix 
\begin{equation*}
A/A[S] = A[\overline{S}]-A[\overline{S},S] A[S]^{-1} A[S,\overline{S}].
\end{equation*} 
Notice that if $A[S]$ is invertible, then $\det(A) = \det(A[S])\det(A/A[S])$.
If $m$ and $n$ are nonnegative integers, we denote the space of all $m\times n$ real matrices by $M_{m,n}$.
If $n$ is a nonnegative integer, we denote the set of all symmetric $n\times n$ real matrices by $S_n$.

\section{Stable Northeast Lemma for signed graphs}\label{sec:inertiaset}

In this section, we show that if $(p,q)\in \Is(G,\Sigma)$ and $p+q < \lvert V(G)\rvert$, then the pair up (to the north) and the pair to the right (to the east) also belong to $\Is(G,\Sigma)$. This theorem was shown for the inertia sets of graphs by Barrett et al. \cite{MR2547901}.

\begin{figure}[h]
\begin{center}
\includegraphics[width=0.7\textwidth]{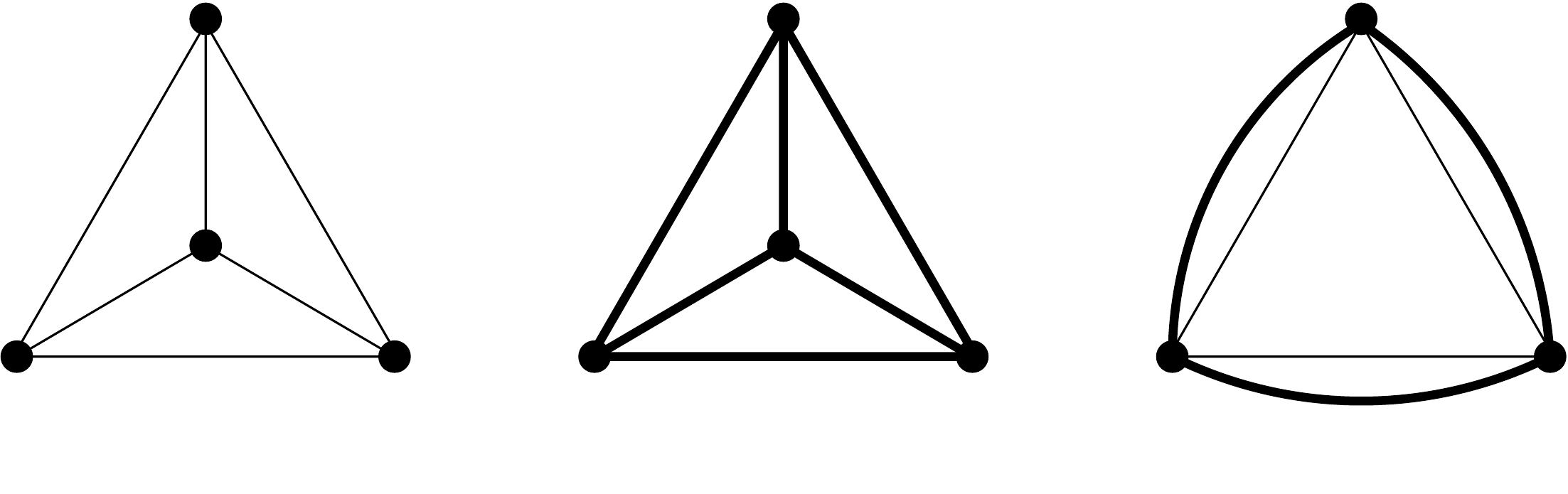}
\end{center}
\caption{$K_4^e$, $K_4^o$, and $K_3^=$}\label{fig:signedgraphs}
\end{figure}

Let us first introduce some signed graphs. By $K_n^e$ and $K_n^o$ we denote the signed graphs $(K_n, \emptyset)$ and $(K_n,E(K_n))$, respectively. By $K_n^=$, we denote the signed graph $(G,\Sigma)$, where $G$ is the graph obtained from $K_n$ by adding to each edge an edge in parallel, and where $\Sigma$ is the set of edges of $K_n$. By $C_n^e$ we denote the signed graph $(C_n,\emptyset)$, and 
by $C_n^o$ we denote the signed graph $(C_n,\{e\})$, where $e$ is an edge of $C_n$.
In Figure~\ref{fig:signedgraphs}, we have depicted $K_4^e$, $K_4^o$, and $K_3^=$. Here a bold edge is an odd edge.
By $K_4^d$, we denote the signed graph $(K_4,\{e\})$, where $e$ is an edge of $K_4$. By $K_{2,3}^e$ and $K_{2,3}^0$, we denote the signed graphs $(K_{2,3},\emptyset)$ and $(K_{2,3},\{e\})$, where $e$ is an edge of $K_{2,3}$, respectively.
  
We have the following observation: if $(p,q) \in \Is(G,\Sigma)$, then $(q,p) \in \Is(G,E(G)\setminus \Sigma)$. A similar statement holds for $\I(G,\Sigma)$. Hence, if $(G,\Sigma)$ is isomorphic to a signed graph that is sign-equivalent to $(G,E(G)\setminus \Sigma)$, then $(q,p)\in \Is(G,\Sigma)$ if $(p,q)\in \Is(G,\Sigma)$, and, of course, a similar statement holds for $\I(G,\Sigma)$. Among the signed graphs that we introduced above, $K_n^=$, $K_{2,3}^e$, $K_{2,3}^o$, $K_4^d$, $C_n^e$, and $C_n^o$, with $n$ even, have this property.

\begin{lemma}[Stable Northeast Lemma for Signed Graphs]\label{lem:snortheast}
Let $(G,\Sigma)$ be a signed graph with $n$ vertices. If $(p,q) \in \Is(G,\Sigma)$ and $p+q<n$, then $(p+1,q),(p,q+1)\in \Is(G,\Sigma)$.
\end{lemma}
\begin{proof}
Let $A\in S(G,\Sigma)$ have the SAP,  $p$ positive and $q$ negative eigenvalues, and
let $k=p+q$. There is a $k\times k$ diagonal matrix $D$ with $p$ positive and $q$ negative numbers on the diagonal, and a $k\times n$ matrix $U$ with $\rank(U)=k$ such that $A = U^T D U$. Denote by $e_i$ the vector whose $i$th coordinate is $1$ and all other coordinates are $0$. Since $k<n$, there is an $1\leq i\leq n$ such that
\begin{equation*}
\begin{bmatrix}
U\\
e_i^T
\end{bmatrix}
\end{equation*}
has rank $k+1$. Let
\begin{equation*}
B = A + e_i^{} e_i^T.
\end{equation*}
Since $A\in S(G,\Sigma)$, clearly $B\in S(G,\Sigma)$.
Since
\begin{equation*}
B = \begin{bmatrix}
U^T & e_i
\end{bmatrix}
\begin{bmatrix}
D & 0\\
0 & 1
\end{bmatrix}
\begin{bmatrix}
U\\
e_i^T
\end{bmatrix},
\end{equation*}
$B$ has, by Sylvester's Law of Inertia, $p+1$ positive and $q$ negative eigenvalues.

It remains to show that $B$ has the SAP\@. The null space of $B$ is equal to the null space
of $
\begin{bmatrix}
U\\
e_i^T
\end{bmatrix}$, so
the null space of $B$ is a subspace of the null space of $A$.
Let $X = [x_{i,j}]$ be a symmetric $n\times n$ matrix with $x_{i,j}=0$ if $i=j$ or $ij\in E$ such that $B X = 0$.
Since each vector in the null space of $B$ is in the null space of $A$, we see that $A X = 0$. As $A$ has the SAP, $X = 0$. Therefore $B$ has the SAP\@.
Thus $(p+1,q)\in \Is(G,\Sigma)$. By taking $B = A-e_i e_i^T$, we see that $(p,q+1)\in \Is(G,\Sigma)$.
\end{proof}

\begin{cor}[Stable Northeast Lemma for Graphs]\label{lem:snortheastgraph}
Let $G$ be a graph with $n$ vertices. If $(p,q) \in \Is(G)$ and $p+q<n$, then $(p+1,q),(p,q+1)\in \Is(G)$.
\end{cor}

The proof of the next lemma is similar to the one for Lemma~\ref{lem:snortheast}.
\begin{lemma}[Northeast Lemma for Signed Graphs]\label{lem:northeast}
Let $(G,\Sigma)$ be a signed graph with $n$ vertices. If $(p,q) \in \I(G,\Sigma)$ and $p+q<n$, then $(p+1,q),(p,q+1)\in \I(G,\Sigma)$.
\end{lemma}

As an illustration of Lemma~\ref{lem:snortheast}, let us determine the inertia sets and the stable inertia sets of the signed graphs $K_n^=$, $K_n^e$, and $K_n^o$. We will determine the inertia sets and stable inertia sets of $C_n^o$, $C_n^e$, $K_{2,3}^e$, $K_{2,3}^o$, and $K_4^d$ in Section~\ref{sec:somegraphs}.

\begin{prop}\label{prop:Kn=}
$\I(K_n^=) = \Is(K_n^=) = \mathbb{N}^2_{[0,n]}$.
\end{prop}
\begin{proof}
The $n\times n$ all-zero matrix belongs to $S(K_n^=)$, has $0$ positive and negative eigenvalues, and has the SAP\@. Hence $(0,0) \in \Is(K_n^=)$. By Lemma~\ref{lem:snortheast}, $\mathbb{N}^2_{[0,n]}\subseteq \Is(K_n^=)$. Since clearly $\I(K_n^=)\subseteq \mathbb{N}^2_{[0,n]}$, we obtain $\I(K_n^=) = \Is(K_n^=) = \mathbb{N}^2_{[0,n]}$.
\end{proof}

To determine the inertia sets of $K_n^e$ and $K_n^o$, we first need a lemma which says that the nullity of a connected bipartite signed graph is at most $1$.
\begin{lemma}\label{lem:smnull1}
Let $(G,\Sigma)$ be a connected bipartite signed graph and let $A\in S(G,\Sigma)$ be positive semidefinite. If $x\in \ker(A)$ is nonzero, then $x$ has only nonzero entries. Especially, $\nullity(A)\leq 1$.
\end{lemma}
\begin{proof}
Since $(G,\Sigma)$ is bipartite, we can resign the signed graph to $(G,\Sigma')$ which has only even edges. Resigning around a vertex $v$ corresponds to multiplying the $v$th row and column of $A$ by $-1$, and to multiplying the $v$ entry of $x$ by $-1$. Hence the resulting matrix $B\in S(G,\Sigma')$ has only nonpositive off-diagonal entries and the resulting vector $y$ belongs to $\ker(B)$. If $a$ sufficiently large, then $-B+a I$ has only nonnegative entries. Notice that $(-B+a I)y = ay$.
Since $a$ is the largest eigenvalue and $-B+a I$ is irreducible, Perron-Frobenius tells us that any nonzero vector in $\ker(-B+aI)$ has all components positive or all components negative. Hence $y$ has only nonzero components, and so has $x$. 
If $\nullity(A)>1$, then there exists a nonzero $x\in \ker(A)$ with at least one component equal to zero. This contradicts that $x$ can have only nonzero components.
\end{proof}

\begin{lemma}\label{lem:possemdef}
If $(G,\Sigma)$ is a connected bipartite signed graph with $n$ vertices, then $(n-2,0)\not\in \I(G,\Sigma)$ and $(n-1,0)\in \Is(G,\Sigma)$.
\end{lemma}
\begin{proof}
By Lemma~\ref{lem:smnull1}, $(n-2,0)\not\in \I(G,\Sigma)$.

To see that $(n-1,0)\in \Is(G,\Sigma)$, let $A\in S(G,\Sigma)$ be positive semidefinite. If the nullity of $A$ is not equal to one, let $\lambda_1$ be the smallest eigenvalue of $A$. Then $A-\lambda_1 I$ is positive semidefinite and has nullity $1$. Hence we may assume that $A$ has nullity $1$. Any matrix $A\in S(G,\Sigma)$ with nullity $1$ has the SAP\@. To see this, let $X=[x_{i,j}]$ be a symmetric matrix satisfying $AX = 0$ and $x_{i,j}=0$ if $i$ and $j$ are adjacent or if $i=j$. Then $X=y y^T$ for some $y\in \ker(A)$ and since $x_{i,i} = 0$ for all vertices $i$, $y = 0$. We can therefore conclude that $(n-1,0)\in \Is(G,\Sigma)$.
\end{proof}

\begin{lemma}\label{lem:deleteS}
Let $(G,\Sigma)$ be a signed graph with $n$ vertices. If $S\subseteq V(G)$ and $(G,\Sigma)\setminus S$ is a connected bipartite signed graph, then $(n-2-|S|,0)\not\in \I(G,\Sigma)$.
\end{lemma}
\begin{proof}
Suppose for a contradiction that $(n-2-\lvert S\rvert,0)\in I(G,\Sigma)$. Let $A\in S(G,\Sigma)$ be a positive semidefinite matrix with nullity $|S|+2$. Then $A[V-S]$ has nullity at least $2$. Since $(G,\Sigma)\setminus S$ is a connected bipartite signed graph, we obtain a contradiction by Lemma~\ref{lem:smnull1}.
\end{proof}

For $A,B\subseteq \mathbb{N}^2$, define
\begin{equation*}
A+B=\{(p_1+p_2,q_1+q_2)\mid (p_1,q_1)\in A, (p_2,q_2)\in B\}.
\end{equation*}

\begin{prop}\label{prop:signedKn}
$\Is(K_1) = \I(K_1) = \mathbb{N}^2_{[0,1]}$. If $n > 1$, then
\begin{equation*}
\Is(K_n^e) = \I(K_n^e) = [\{(0,1)\}+\mathbb{N}^2_{[0,n-1]}]\cup \{(n-1,0),(n,0)\}
\end{equation*}
and 
\begin{equation*}
\Is(K_n^o) = \I(K_n^o) = [\{(1,0)\}+\mathbb{N}^2_{[0,n-1]}]\cup \{(0,n-1),(0,n)\}.
\end{equation*}
\end{prop}
\begin{proof}
The $n\times n$ all $-1$ matrix has one negative, no positive eigenvalues, and has the SAP\@. By the Stable Northeast Lemma, $\{(0,1)\}+\mathbb{N}^2_{[0,n-1]}\subseteq \Is(K_n^e)$.
Since $K_n^e$ is a connected bipartite signed graph, $(n-2,0)\not\in \I(K_n^e)$. Since $\Is(K_n^e)\subseteq \I(K_n^e)$, 
$\Is(K_n^e) = \I(K_n^e) = [\{(0,1)\}+\mathbb{N}^2_{[0,n-1]}]\cup \{(n-1,0),(n,0)\}$. The (stable) inertia set of $K_n^o$ follows immediately from the (stable) inertia set of $K_n^e$.
\end{proof}

If $0\leq k\leq n$, we define $\mathbb{N}^2_{[k,n]} = \{(p,q)\in \mathbb{N}^2\mid k\leq p+q\leq n\}$.

In \cite{MR2547901}, Barrett, Hall, and Loewy showed that $\I(K_n) = \mathbb{N}^2_{[1,n]}$. From Proposition~\ref{prop:signedKn}, we obtain the following extension of their result.

\begin{cor}\label{cor:IsKn}
$\Is(K_n) = \I(K_n) = \mathbb{N}^2_{[1,n]}$.
\end{cor}
\begin{proof}
Since $\Is(K_n^e)\cup \Is(K_n^o) \subseteq \Is(K_n)$ and $\mathbb{N}^2_{[1,n]}\subseteq \Is(K_n^e)\cup \Is(K_n^o)$, $\mathbb{N}^2_{[1,n]}\subseteq \Is(K_n)$. Since $\mr(K_n)=1$, $\I(K_n)\subseteq \mathbb{N}^2_{[1,n]}$. Hence $\Is(K_n) = \I(K_n) = \mathbb{N}^2_{[1,n]}$.
\end{proof}

\section{The Strong Arnold Hypothesis}\label{sec:SAH}

Let $G=(V,E)$ be a graph and $A\in S(G)$. Suppose $R\subseteq V$ such that $A[R]$ is invertible and $A/A[R] = 0$. Let $U$ be an open neighborhood around $A$ in $S(G)$ such that $B[R]$ is invertible for each $B\in U$. Define 
\begin{equation*}
\psi_S : U\to \mathcal{S}_{k}
\end{equation*} 
by 
\begin{equation*}
\psi_R(B) = B/B[R],\quad B\in U.
\end{equation*} 
If it is clear from the context which set $R$ we use, we will write $\psi$ instead of $\psi_R$. 
In this section we show that $A$ has the SAP if and only if the derivative, $D\psi_R(A)$, of $\psi_R$ at $A$ is surjective. First we give an alternative criterion for a matrix $A\in S(G)$ to have the SAP; this is the criterion that Colin de Verdi\`ere used in his papers \cite{CdeV1, CdeV2, CdeV3}. 

If $G$ is a graph with vertex set $\{1,\ldots,n\}$, we denote by 
$T(G)$ the set of all symmetric $n\times n$ matrices $B=[b_{i,j}]$ with $b_{i,j} = 0$ if $i\not=j$ and $i$ and $j$ are not adjacent. 
In the following proof we will use that $\Tr(C D) = \Tr(D C)$  for any $n\times m$ matrix $C$ and any $m\times n$ matrix $D$.

\begin{thm} \label{crit:colindeverdiere}
Let $G$ be a graph with $n$ vertices. A matrix $A\in S(G)$ has the SAP if and only if for each
real symmetric $n\times n$ matrix $N$, there is a matrix $B \in T(G)$ such that $x^T N x = x^T B x$ for all $x \in \ker(A)$.
\end{thm}
\begin{proof}
Let $y_1,\ldots,y_k$ form a basis of $\ker(A)$ and let $Y = \bigl[ y_1 \dots y_k \bigr]$. The statement that for each real symmetric $n\times n$ matrix $N$, there is a matrix $B\in T(G)$ such that $x^T N x = x^T B x$ for each $x\in \ker(A)$ is the same as the statement that for each real symmetric $n\times n$ matrix $N$, there is a matrix $B\in T(G)$ such that $Y^T N Y = Y^T B Y$.

Suppose that $A$ does not have the SAP\@. Then there exists a nonzero real symmetric matrix $X$ with $x_{i,i} = 0$ and $x_{i,j} = 0$ if $i$ and $j$ are adjacent in $G$, such that $A X =0$.
Since $x_{i,j} = 0$ if $i$ and $j$ are adjacent or if $i=j$, $\Tr(X K) = 0$ for all $K\in T(G)$. Since each column of $X$ belongs to $\ker(A)$, there exists a nonzero symmetric $k\times k$ matrix $C$ such that $X = Y C Y^T$. Hence $\Tr(Y C Y^T K) = 0$ for all $K\in T(G)$, and therefore 
\begin{equation*}
\Tr(C Y^T K Y) = 0\quad\text{for all $K\in T(G)$.}
\end{equation*}
Let $N$ be a symmetric $n\times n$ matrix such that $Y^T N Y = C$. If there would exist a symmetric matrix $B\in T(G)$ such that $C = Y^T B Y$, then $\Tr(Y^T B Y Y^T K Y)=0$ for all $K\in T(G)$.
Especially, $\Tr(Y^T B Y Y^T B Y) = 0$, and hence $C = Y^T B Y = 0$. This contradiction shows that there is no such $B\in T(G)$.

Conversely, suppose that there exists a real symmetric $n\times n$ matrix $N$ such that there is no matrix $B\in T(G)$ for which $x^T N x = x^T B x$ for each $x\in \ker(A)$. 
Then the space of all matrices $Y^T B Y$ with $B \in T(G)$ is not equal to the space of all real symmetric $k\times k$ matrices. This means that there exists a real symmetric $k\times k$ matrix $C$ such that $\Tr(Y^T B Y C) = 0$ for all $B \in T(G)$.
Since $\Tr(Y^T B Y C) = \Tr(B Y C Y^T) = 0$, we see that
$Y C Y^T$ is orthogonal to each matrix in $T(G)$. So
$X = [x_{i,j}] = Y C Y^T$ is a nonzero matrix with $x_{i,j} = 0$ if $i=j$ or if $i$ and $j$ are adjacent in $G$ for
which $A X =0$.
\end{proof}

 Let $L$ be a subspace of $M_{m,n}$. If $f : L\times L\to M_{m,n}$ is defined by $f(A,B) = A+B$, then the derivative of $f$ at $(A,B)$ is the linear map $Df(A,B) : L\times L\to M_{m,n}$ defined by
\begin{equation}\label{eq:dersum}
Df(A,B)(C,D) = C+D.
\end{equation} 
Let $L_1$ and $L_2$ be subspaces of $M_{m,n}$ and $M_{n,p}$, respectively. 
If $g :  L_1\times L_2\to M_{m,p}$ is defined by $g(A,B) = AB$, then the derivative of $g$ at $(A,B)$ is the linear map $Dg(A,B) : L_1\times L_2\to M_{m,p}$ defined by
\begin{equation}\label{eq:derprod}
Dg(A,B)(C,D) = CB + AD.
\end{equation} 
Let $L$ be a subspace of the space of all $n\times n$ matrices and let $U$ be the open set of all invertible matrices in $L$.
If $h : U\to M_n$ is defined by $h(A) = A^{-1}$, then 
the derivative of $h$ at $A$ is the linear map $Dh(A): L\to M_n$ defined by
\begin{equation}\label{eq:derinv}
Dh(A)(C) = -A^{-1}CA^{-1}.
\end{equation} 

Define 
\begin{equation*}
\gamma_R(A) = \begin{bmatrix}
-A[R]^{-1}A[R,\overline{R}]\\
I_{k}
\end{bmatrix}.
\end{equation*}
The columns of $\gamma_R(A)$ form a basis of $\ker(A)$.
The mapping $\psi_R$ is $C^{1}$ and, by Equations~\ref{eq:dersum},~\ref{eq:derprod}, and ~\ref{eq:derinv}, the derivative at $A$ is the linear map 
\begin{equation*}
D\psi_R(A) : T(G)\to \mathcal{S}_{k}
\end{equation*}
defined by 
\begin{equation*}
D\psi_R(A)(C) = \gamma_R(A)^T C \gamma_R(A).
\end{equation*}
(A $C^1$ map is a differentiable map whose derivative is continuous.)

Finally, we present our main result of this section; it is a new result that will be important in the next section.

\begin{thm}\label{thm:SAHsubmersion}
Let $A \in S(G)$ and let $R\subseteq V(G)$ such that $A[R]$ is invertible and $\psi_R(A) = 0$. Then $A$ has the SAP if and only if $D\psi_R(A)$ is surjective.
\end{thm}
\begin{proof}
Suppose $A\in S(G)$ has the SAP\@. Let $B\in \mathcal{S}_{k}$. The symmetric $n\times n$ matrix
\begin{equation*}
B' = \begin{bmatrix}
0 & 0\\
0 & B
\end{bmatrix}
\end{equation*}
satisfies $\gamma_R(A)^T B' \gamma_R(A) = B$. Since $A$ has the SAP, there exists a matrix $C\in T(G)$ such that $\gamma_R(A)^T C \gamma_R(A) = \gamma_R(A)^T B' \gamma_R(A)$. Thus $\gamma_R(A)^T C\gamma_R(A) = B$. 

Conversely, suppose $D\psi_R(A)$ is surjective. Let $B'$ be a symmetric $n\times n$ matrix and let $B = \gamma_R(A)^T B' \gamma_R(A)$. There exists a matrix $C\in T(G)$ such that $D\psi_R(A)(C) = \gamma_R(A)^T C\gamma_R(A) = B$. Thus $\gamma_R(A)^T B' \gamma_R(A) = \gamma_R(A)^T C \gamma_R(A)$.
\end{proof}

Let $A\in S(G,\Sigma)$ and let $R\subseteq V(G)$ such that $A[R]$ is invertible and $\psi_R(A) = 0$. Suppose that $D\psi_R(A)$ is surjective. Since $D\psi_R(A)$ is a continuous function in $A$, there exists a neighborhood $U$ of $A$ in $S(G,\Sigma)$ such that $D\psi_R(B)$ is surjective for all $B\in U$.

\section{The stable inertia set and subgraphs and contracting edges}\label{sec:minorsstableinertiaset}

This is the main section of the paper. If $k \in \mathbb{N}$, then by $\mathbb{N}^2_k$ we denote the set $\{(p,q)\in \mathbb{N}^2\mid p+q=k\}$. Let $(G,\Sigma)$ be a signed graph with $n$ vertices and let $G'$ be a subgraph of $G$ with $m$ vertices. In this section we prove that if $(H,\Omega)$ is obtained from $(G',\Sigma\cap E(G'))$ by contracting $s$ even and $t$ odd edges, then 
\begin{equation*}
\Is(H,\Omega)+\mathbb{N}^2_{n-m} + \{(s,t)\}\subseteq \Is(G,\Sigma).
\end{equation*}
First we state the Implicit Function Theorem and give some lemmas.

Let $(x_0,y_0) \in \mathbb{R}^n\times\mathbb{R}^m$ and $W\subseteq \mathbb{R}^n\times\mathbb{R}^m$ be a neighborhood of $(x_0,y_0)$. If $F : W\to \mathbb{R}^m$ is differentiable, then $D_1 F(x_0,y_0)$ denotes the derivative of $j(x) = F(x,y_0)$ at $x = x_0$ and $D_2 F(x_0,y_0)$ the derivative of $k(y) = F(x_0,y)$ at $y=y_0$.

\begin{thm}[Implicit Function Theorem]\cite{gast.26511419740101}
Let $W\subseteq \mathbb{R}^n\times \mathbb{R}^m$ be an open set and $(x_0,y_0) \in W$. Let $F : W\to \mathbb{R}^m$ be a $C^1$ map and $c = F(x_0,y_0)$. Suppose $(x_0,y_0)\in W$ is such that the linear operator $D_2 F(x_0,y_0) : \mathbb{R}^m\to \mathbb{R}^m$ is invertible. Then there are open sets $U\subseteq \mathbb{R}^n$ and $V\subseteq \mathbb{R}^m$ with 
\begin{equation*}
(x_0,y_0)\in U\times V\subseteq W
\end{equation*} 
and a unique $C^1$ map $g : U\to V$ such that 
\begin{equation*}
F(x,g(x)) = c
\end{equation*} 
for all $x\in U$, and moreover, $F(x,y) \not= c$ if $(x,y)\in U\times V$ and $y\not=g(x)$.
\end{thm}

We apply the Implicit Function Theorem in the proof of the following lemma.

\begin{lemma}\label{lem:neighborSAH}
Let $W\subseteq \mathbb{R}\times\mathbb{R}^n$ be an open set and $(0,a) \in W$. Let $F : W\to \mathbb{R}^m$ be a $C^1$ map, and define $f(x) = F(0,x)$. If $f(a) = 0$ and $Df(a)$ is surjective,  then there exists a $\delta > 0$ such that for each $h\in \mathbb{R}$ with $|h|<\delta$, there exists a $(h,b)\in W$ such that $F(h,b) = 0$ and $DF(h,b)$ is surjective.
\end{lemma}

\begin{proof}
Since $Df(a)$ is surjective, the Jacobian matrix $Jf(a)$ of $f$ at $a$ has rank $m$. We may assume that the last $m$ columns of $Jf(a)$ span an invertible matrix. We write $\mathbb{R}^n = \mathbb{R}^{n-m}\times\mathbb{R}^m$ and $a = (a_1,a_2)$, where $a_1\in \mathbb{R}^{n-m}$ and $a_2\in \mathbb{R}^m$.
Then $D_2 f(a_1,a_2)$ is an invertible linear operator, which implies that $D F(0,a)$ is surjective. By the Implicit Function Theorem, there are open sets $U\subseteq \mathbb{R}^{n-m+1}$ and $V\subseteq \mathbb{R}^m$ such that $(0,a_1) \in U$, $a_2\in V$ and $U\times V\subseteq W$ and a unique $C^1$ map $G : U\to V$ such that $F(x,G(x)) = 0$ for all $x\in U$ and $F(x,y)\not=0$ if $(x,y)\in U\times V$ and $y\not=G(x)$. Let $\delta>0$ for which there exists a neighborhood $U'$ of $a_1$ such that $(-\delta,\delta)\times U'\subseteq U$. If $h\in (-\delta,\delta)$, choose $b_1\in U'$ and let $b = (b_1,G(h,b_1))$. Then $F(h,b) = 0$.
\end{proof}

The following lemma is taken from \cite{Holst2007b}; for convenience we include a proof.
\begin{lemma}\label{lem:schurreduction}
Let $A=[a_{i,j}]$ be a nonzero symmetric $n\times n$ matrix. Then there is an $S\subseteq \{1,\ldots,n\}$ such that $A[S]$ is invertible and $A/A[S] = 0$.
\end{lemma}
\begin{proof}
Take $S \subseteq \{1,\ldots,n\}$ such that $A[S]$ is invertible and
$\lvert S\rvert$ is as large as possible. (It is clear that such an $S$ exists.) Let $B = [b_{i,j}] = A/A[S]$.
If $b_{i,i} \not= 0$ for some $i \not\in S$, then
\begin{equation*}
\begin{aligned}
\det(A[S\cup \{i\}]) & = \det(A[S]) \det(A[S\cup\{i\}]/A[S]) \\
        & = \det(A[S]) b_{i,i} \\
& \not= 0.
\end{aligned}
\end{equation*}
This contradicts the maximality of $\lvert S\rvert$. Hence $b_{i,i} = 0$ for each $i \not\in S$.

If there are vertices $i,j \not\in S$ such that
$b_{i,j}\not=0$, then the principal submatrix
\begin{equation*}
K = \begin{bmatrix}
0 & b_{i,j} \\ b_{j,i} & 0
\end{bmatrix}
\end{equation*}
of $B$ has $\det(K)\not=0$.
So
\begin{equation*}
\begin{aligned}
\det(A[S\cup \{i,j\}]) & =  \det(A[S]) \det(A[S\cup\{i,j\}]/A[S])\\
        & =  \det(A[S]) \det(K) \\
        & \not= 0,
\end{aligned}
\end{equation*}
which contradicts the maximality of $\lvert S\rvert$.
Thus  $B = 0$.
\end{proof}

\begin{thm} \label{thm:stableinertiasubgraph}
Let $(G,\Sigma)$ be a signed graph with $n$ vertices and let $(H, \Sigma\cap E(H))$ be a signed subgraph of $(G,\Sigma)$ with $m$ vertices. Then
\begin{equation*}
\Is(H,\Sigma\cap E(H)) +\mathbb{N}^2_{n-m}\subseteq \Is(G,\Sigma).
\end{equation*}
\end{thm}

\begin{proof}
We may assume that $(G,\Sigma)$ has vertex set $\{1,\dots,n\}$.
As each signed subgraph $(H,\Sigma\cap E(H))$ of $(G,\Sigma)$ arises by a series of deletions of edges and isolated vertices, we may assume that either $H = G\setminus\{v\}$
for an isolated vertex $v$ of $G$ or $H=G\setminus e$ for an edge $e$.

Assume $(H,\Sigma)$ arises from $(G,\Sigma)$ by deleting an isolated vertex $v$; we assume that $v=1$.
Let $A = [a_{i,j}] \in S(H,\Sigma)$ have the SAP and have $p$ positive and $q$ negative eigenvalues. The matrix
\begin{equation*}
M = \begin{bmatrix}
1 & 0\\
0 & A
\end{bmatrix}
\end{equation*}
belongs to $S(G,\Sigma)$ and has $p+1$ positive and $q$ negative eigenvalues. We claim that $M$ has the SAP\@. Each vector $x\in \ker(M)$ has $x_1 = 0$ and $x[V(H)]\in \ker(A)$. Let $C$ be an arbitrary real symmetric $n\times n$ matrix,  $C' = C[V(H)]$, and $y = x[V(H)]$. As $A$ has the SAP, there exists a matrix $B' = [b'_{i,j}] \in T(H)$ such 
$y^T C' y = y^T B' y$. The matrix
\begin{equation*}
B = \begin{bmatrix}
0 & 0\\
0 & B'
\end{bmatrix}
\end{equation*}
belongs to $T(G)$. 
We get $x^T C x = y^T C' y = y^T B' y = x^T B x$ for all $x \in \ker(M)$, showing that $M$ has the SAP\@.

Assume next that $(H,\Sigma\cap E(H))$ arises from $(G,\Sigma)$ by deleting an even edge $e = uw$;  the case where $e$ is odd is similar, except that in that case we have to take $h>0$.
We may assume that $u=1$ and $w=2$. Let $R := \{3,\ldots,n\}$. Define a function 
\begin{equation*}
f: \mathbb{R} \times S(H,\Sigma\cap E(H))\to \mathcal{S}_n
\end{equation*}
by
\begin{equation*}
f(h,K) := \begin{bmatrix}
k_{1,1} & h & K[1, R]\\
h & k_{2,2} & K[2, R]\\
K[R,1] & K[R, 2] & K[R]
\end{bmatrix},
\end{equation*}
where $K = [k_{i,j}] \in S(H,\Sigma\cap E(H))$, and let $f_h(K) = f(h,K)$. 
So $f_0$ is the identity on $S(H,\Sigma\cap E(H))$. 

Let $A = [a_{i,j}] \in S(H,\Sigma\cap E(H))$ have the SAP and have $p$ positive and $q$ negative eigenvalues. By Lemma~\ref{lem:schurreduction}, there exists a subset $S\subseteq V$ for which $A[S]$ is invertible and $\psi_S(A) = 0$. Since $A$ has the SAP, $D\psi_S(A)$ is surjective, by Theorem~\ref{thm:SAHsubmersion}. Since eigenvalues continuously depend on the entries of the matrix, there exists an open set $U$ of $\mathbb{R}\times S(H,\Sigma\cap E(H))$ containing $(0,A)$ such that for all $(h,B)\in U$, $f(h,B)[S]$ has the same inertia as $A[S] = f(0,A)[S]$. 
Define $g : U\to \mathcal{S}_k$ by 
$g = \psi_S\circ f_{|U}$
and define $g_0(x) = g(0,x)$. Then $g_0(A) = 0$ and $D g_0(A)$ is surjective. By Lemma~\ref{lem:neighborSAH}, there exists a $\delta>0$ such that for all $h_1$ with $\lvert h_1\rvert<\delta$, there exists a $C_1\in S(H,\Sigma\cap E(H))$ such that $(h_1,C_1)\in U$, $g(h_1,C_1)=0$, and $D g(h_1,C_1)$ is surjective. Let $\delta <  h < 0$ and $C \in S(H,\Sigma\cap E(H))$ such that $(h,C)\in U$, $g(h,C)=0$ and $D g(h,C)$ is surjective.
Let $M = f(h,C)$. Then $M\in S(G,\Sigma)$. Since $D g(h,C)$ is surjective, $D \psi_S(M)$ is surjective. Hence $M$ has the SAP\@. Since $(h,C) \in U$, $M[S]$ has the same inertia as $A[S]$. Thus $\Is(H,\Sigma\cap E(H)) \subseteq \Is(G,\Sigma)$.
\end{proof}  

\begin{lemma}\label{lem:twoband}
If $(G,\Sigma)$ is a signed graph with $n>0$ vertices, then $\mathbb{N}^2_{[n-1,n]}\subseteq \Is(G,\Sigma)\subseteq \I(G,\Sigma)$.
\end{lemma}
\begin{proof}
Since $\Is(K_1) = \mathbb{N}^2_{[0,1]}$ and by deleting $n-1$ vertices from $(G,\Sigma)$ we obtain a single vertex, $\mathbb{N}^2_{[0,1]}+\mathbb{N}^2_{n-1} = \mathbb{N}^2_{[n-1,n]}\subseteq \Is(G,\Sigma)$.
\end{proof}

Barrett et al. \cite{MR2547901} showed that for any graph with $n$ vertices, $\mathbb{N}^2_{[n-1,n]}\subseteq \I(G)$. Lemma~\ref{lem:twoband} is a strengthening of their result, as for any signed graph $(G,\Sigma)$, $\mathbb{N}^2_{[n-1,n]}\subseteq \Is(G,\Sigma)\subseteq \Is(G)\subseteq \I(G)$.

If $G$ is a graph without loops and $e$ is an edge of $G$, then $G/e$ denotes the graph obtained from $G$ by contracting $e$.

\begin{thm} \label{thm:stableinertiaminor}
Let $(G,\Sigma)$ be a signed graph and $e\in E(G)$.
If $e$ is odd, then 
\begin{equation*}
\Is(G/e,\Sigma\setminus\{e\})+\{(0,1)\}\subseteq \Is(G,\Sigma).
\end{equation*}
If $e$ is even, then 
\begin{equation*}
\Is(G/e,\Sigma)+\{(1,0)\}\subseteq \Is(G,\Sigma).
\end{equation*}
\end{thm}

\begin{proof}
We may assume that $(G,\Sigma)$ has vertex set $\{1,\dots,n\}$.

We give the proof only for the case that $e=uw$ is an even edge; the case that $e$ is odd is similar, except that $h<0$ when $e$ is odd. We may assume that $u=n-1$ and $w=n$. Let $v$ be the new vertex of $G/e$; so $v=n-1$.
  
Let $R=\{1,\ldots,n-2\}$. Let $\mathcal{Z}$ be the subset of all matrices 
$K=[k_{i,j}] \in S(G,\Sigma)$ with $k_{n,n}=0= k_{n,n-1}$. 
Define a function 
\begin{equation}
f : \oR\times \mathcal{Z} \longrightarrow \mathcal{S}_{n-1},
\end{equation}
by
\begin{equation}
f(h,K) = \begin{bmatrix}
K[R] - h K[R, n] K[n, R] & K[R, n] + K[R, n-1]\\
K[n,R] + K[n-1,R] & k_{n-1,n-1}
\end{bmatrix},
\end{equation}
and let $f_h(K) = f(h,K)$. 

Let $A=[a_{i,j}]\in S(G/e,\Sigma)$ have the SAP and have $p$ positive and $q$ negative eigenvalues. By Lemma~\ref{lem:schurreduction}, there exists a subset $S\subseteq \{1,\ldots,n-1\}$ for which $A[S]$ is invertible and $\psi_S(A) = 0$. Since $A$ has the SAP, $D \psi_S(A)$ is surjective by Theorem~\ref{thm:SAHsubmersion}. The image of $f_0$ is a superset of $S(G/e,\Sigma)$, so there is a $P \in \oR\times \mathcal{Z}$ such that $f_0(P) = A$. Since eigenvalues continuously depend on the entries of the matrix, there exists an open set $U$ of $\mathbb{R}\times\mathcal{Z}$ containing $(0,P)$ such that for all $(h,B)\in U$, $f(h,B)[S]$ has the same inertia as $A[S] = f(0,P)[S]$. Define $g : U\to \mathcal{S}_k$ by $g = \psi_S\circ f_{|U}$ and define $g_0(x) = g(0,x)$.  Then $g_0(P) = 0$.  Since $T(G/e)$ is a subset of the image of $Df_{|U}$, $D g_0(P)$ is surjective. By Lemma~\ref{lem:neighborSAH}, there exists a $\delta>0$ such that for all $h_1$ with $\lvert h_1\rvert<\delta$, there exists a $C_1\in S(G/e,\Sigma)$ such that $(h_1,C_1)\in U$, $g(h_1,C_1)=0$, and $D g(h_1,C_1)$ is surjective. Let $0 < h<\delta$ and $C\in S(G/e,\Sigma)$ such that $(h,C)\in U$, $g(h,C)=0$, and $D g(h,C)$ is surjective. Let
\begin{equation*}
M = \begin{bmatrix}
C[R] & C[R, n-1] & C[R, n]\\
C[n-1,R] & c_{n-1,n-1}+\frac{1}{h} & -\frac{1}{h}\\
C[n, R] & -\frac{1}{h} & \frac{1}{h} 
\end{bmatrix}.
\end{equation*}
Then $\psi_{\{n\}}(M) = f(h,C)$, so $g(h,C) = \psi_S(\psi_{\{n\}}(M))$. Since $(h,C)\in U$, $f(h,C)[S]$ has the same inertia as $A[S] = f(0,P)[S]$. Hence $M$ has $p+1$ positive and $q$ negative eigenvalues. Since $D g(h,C)$ is surjective, $D (\psi_S\circ\psi_{\{n\}})(M) = D\psi_{S\cup\{n\}}(M)$ is surjective. Hence $M$ has the SAP\@. Thus $\Is(G/e,\Sigma)+\{(1,0)\}\subseteq \Is(G,\Sigma)$.
\end{proof}  

From the previous theorem we obtain the following theorem.

\begin{thm}\label{thm:contractingedge}
Let $(G,\Sigma)$ be a signed graph, let $e$ be an even edge, and let $u$ be an end of $e$. Then
$[\{(1,0)\}+\Is(G/e,\Sigma)]\cup [\{0,1\}+\Is(G/e,\Sigma\Delta (\delta(u)\setminus\{e\})]\subseteq \Is(G,\Sigma)$.
\end{thm}

Recall that a signed graph $(H,\Sigma')$ is a \emph{minor} of a signed graph $(G,\Sigma)$ if $(H,\Sigma')$ is sign-equivalent to a signed graph that can be obtained by contracting a sequence of even edges in a subgraph of $(G,\Sigma)$.
\begin{cor}
If $(G',\Sigma')$ is a minor of $(G,\Sigma)$, then $\nu(G',\Sigma')\leq \nu(G,\Sigma)$.
\end{cor}
\begin{proof}
Let $(G,\Sigma)$ have $n$ vertices. As resigning on a subset $U\subseteq V(G)$ does not change $\nu(G,\Sigma)$, it suffices to prove this for the cases that $(G',\Sigma')$ is obtained from $(G,\Sigma)$ by deleting an edge or an isolated vertex, or contracting an even edge.

Suppose that $(G',\Sigma')$ is obtain from $(G,\Sigma)$ by deleting an isolated vertex. Then $(G',\Sigma')$ has $n-1$ vertices. Since $(n-1-\nu(G',\Sigma'))\in \Is(G',\Sigma)$, we obtain by Theorem~\ref{thm:stableinertiasubgraph} that $\{(n-1-\nu(G',\Sigma'),0)\}+\mathbb{N}^2_1\subseteq \Is(G,\Sigma)$. Hence $(n-\nu(G',\Sigma'),0)\in \Is(G,\Sigma)$ and so $\nu(G',\Sigma')\leq \nu(G,\Sigma)$.

Suppose next that $(G',\Sigma')$ is obtained from $(G,\Sigma)$ by deleting an edge. Then $(G',\Sigma')$ has $n$ vertices. Since $(n-\nu(G',\Sigma'),0)\in \Is(G',\Sigma')$, we obtain by Theorem~\ref{thm:stableinertiasubgraph} that $(n-\nu(G',\Sigma'),0)\in \Is(G,\Sigma)$. Hence $\nu(G',\Sigma')\leq \nu(G,\Sigma)$.

Finally, suppose that $(G',\Sigma')$ is obtained from $(G,\Sigma)$ by contracting an even edge. Then $(G',\Sigma')$ has $n-1$ vertices. Since $(n-1-\nu(G',\Sigma'),0)\in \Is(G',\Sigma')$, we obtain by Theorem~\ref{thm:stableinertiaminor}  that $\{(n-1-\nu(G',\Sigma'),0)\}+\{(1,0)\}\in \Is(G,\Sigma)$. Hence $(n-\nu(G',\Sigma'),0)\in \Is(G,\Sigma)$ and so $\nu(G',\Sigma')\leq \nu(G,\Sigma)$.
\end{proof}

In the same way one can prove the following corollary.
\begin{cor}
If $H=(W,F)$ is a minor of $G=(V,E)$ and $(G,\Sigma)$ is a signed graph, then $\xi(H, \Sigma\cap F)\leq \xi(G,\Sigma)$.
\end{cor}

\begin{cor}\cite{CdeV2}
If $H$ is a minor of $G$, then $\mu(H)\leq \mu(G)$.
\end{cor}
\begin{proof}
It suffices to prove this for the cases that $H$ is obtained from $G$ by deleting an edge or an isolated vertex, or contracting an edge. Let $G$ have $n$ vertices.

Suppose that $H$ is obtained from $G$ by deleting an isolated vertex. Then $H$ has $n-1$ vertices. Since $(n-1-\mu(H)-1,1)\in \Is(H,\emptyset)$, we obtain by Theorem~\ref{thm:stableinertiasubgraph} that $\{(n-2-\mu(H),1)\}+\mathbb{N}^2_1\subseteq \Is(G,\emptyset)$. Hence $(n-1-\mu(H),1)\in \Is(G,\emptyset)$ and so $\mu(H)\leq \mu(G)$.

Suppose next that $H$ is obtained from $G$ by deleting an edge. Then $H$ has $n$ vertices. Since $(n-\mu(H)-1,1)\in \Is(H,\emptyset)$, we obtain by Theorem~\ref{thm:stableinertiasubgraph} that $(n-1-\mu(H),1)\in \Is(G,\emptyset)$. Hence $\mu(H)\leq \mu(G)$.

Finally, suppose that $H$ is obtained from $G$ by contracting an edge. Then $H$ has $n-1$ vertices. Since $(n-1-\mu(H)-1,1)\in \Is(H,\emptyset)$, we obtain by Theorem~\ref{thm:stableinertiaminor}  that $\{(n-2-\mu(H),1)\}+\{(1,0)\}\in \Is(G,\emptyset)$. Hence $(n-1-\mu(H),1)\in \Is(G,\emptyset)$ and so $\mu(H)\leq \mu(G)$.
\end{proof}

\begin{cor}\label{cor:graphminormonotonicity}
Let $H$ be a minor of $G$ and let $k=|G|-|H|$. Then $\Is(H)+\mathbb{N}^2_k\subseteq \Is(G)$.
\end{cor}
\begin{proof}
It suffices to prove this for the cases that $H$ is a subgraph of $G$ or that $H$ is obtained from $G$ by contracting an edge.  Let $(p,q)\in \Is(H)$.
Then there exists a subset $\Sigma\subseteq E(H)$ such that $(p,q)\in \Is(H,\Sigma)$. If $H$ is a subgraph of $G$, then, by Theorem~\ref{thm:stableinertiasubgraph}, $\{(p,q)\}+\mathbb{N}^2_k\subseteq \Is(G,\Sigma)\subseteq \Is(G)$. 
If $H$ is obtained from $G$ by contracting an edge, then, by Theorem~\ref{thm:stableinertiaminor}, $\{(p,q)\} + \mathbb{N}^2_1\subseteq \Is(G,\Sigma)\subseteq \Is(G)$. 
\end{proof}

\begin{cor}\cite{CdeV3}
If $H$ is a minor of $G$, then $\nu(H)\leq \nu(G)$.
\end{cor}
\begin{proof}
There exists a subset $\Sigma\subseteq E(H)$ such that $\nu(H) = \nu(H,\Sigma)$. Since $(H,\Sigma)$ is a minor of $(G,\Sigma)$, $\nu(H,\Sigma)\leq \nu(G,\Sigma)$. As $\nu(G,\Sigma)\leq\nu(G)$, we obtain that $\nu(H)\leq \nu(G)$.
\end{proof}

In the same way one can prove the following corollary.
\begin{cor}\cite{BarFalHog2005a}
If $H$ is a minor of $G$, then $\xi(H)\leq \xi(G)$.
\end{cor}


\section{The stable inertia sets of some other graphs and signed graphs}\label{sec:somegraphs}

In this section we determine the inertia sets of some more graphs and signed graphs. We first determine the stable inertia sets of trees. Since for any tree $T$ and any $\Sigma\subseteq E(T)$, the signed graph $(T,\Sigma)$ is sign-equivalent to $(T,\emptyset)$, these results immediately give the stable inertia set of signed graphs $(T,\Sigma)$, where $T$ is a tree and $\Sigma\subseteq E(T)$. 

\begin{prop}
If $n\geq 2$, then $\Is(P_n) = \I(P_n) = \oN^2_{[n-1,n]}$, where $P_n$ is a path on $n$ vertices.
\end{prop}
\begin{proof}
By Lemma~\ref{lem:twoband}, $\oN^2_{[n-1,n]}\subseteq \Is(P_n)\subseteq \I(P_n)$. As $\mr(P_n) = n-1$ (see \cite{FalHog2007,MR0244285}), $\I(P_n) \subseteq \oN^2_{[n-1,n]}$. Hence $\I(P_n) = \Is(P_n) = \oN^2_{[n-1,n]}$.
\end{proof}

\begin{prop}\label{prop:K13}
$\Is(K_{1,3}) = \I(K_{1,3}) = \oN^2_{[3,4]} \cup \{(1,1)\}$.
\end{prop}
\begin{proof}
By Lemma~\ref{lem:twoband}, $\oN^2_{[3,4]}\subseteq \Is(K_{1,3})\subseteq \I(K_{1,3})$.
Since $\mu(K_{1,3}) = 2$, $(1,1) \in \Is(K_{1,3},\emptyset)$, and hence $(1,1)\in \Is(K_{1,3})$.
Since $\mr(K_{1,3}) = 2$, $\I(K_{1,3})\subseteq \oN^2_{[2,4]}$. Since $\mr_+(K_{1,3}) = 3$,
$\I(K_{1,3}) \subseteq \oN^2_{[3,4]} \cup \{(1,1)\}$. Hence $\Is(K_{1,3}) = \I(K_{1,3}) = \oN^2_{[3,4]} \cup \{(1,1)\}$.
\end{proof}

\begin{prop}\label{prop:treestabinertia}
If $T$ is a tree with $n$ vertices which has at least one vertex of degree $\geq 3$, then $\Is(T) = (\oN^2_{[3,4]} \cup \{(1,1)\}) + \oN^2_{n-4}$.
\end{prop}
\begin{proof}
Since $K_{1,3}$ is a minor of $T$, $(\oN^2_{[3,4]} \cup \{(1,1)\}) + \oN^2_{n-4} \subseteq\Is(T)$ by Theorem~\ref{cor:graphminormonotonicity}. Since $\xi(T) \leq 2$ (see \cite{BarFalHog2005a}), $\Is(T) \subseteq \oN^2_{[n-2,n]}$. Since $M_+(T) = 1$, $\I_s(T) = (\oN^2_{[3,4]} \cup \{(1,1)\}) + \oN^2_{n-4}$.
\end{proof}

In general, $\I(T)\not=\Is(T)$ for a tree $T$. For example, if $T=K_{1,4}$, then $M(K_{1,4}) = 3$, while $\xi(K_{1,4}) = 2$, so $\Is(K_{1,4})$ is a strict subset of $\I(K_{1,4})$.

\begin{thm}
Let $(G,\Sigma)$ be a signed graph and let $v$ be a vertex of degree one. Let $w$ be the vertex adjacent to $v$. Let $(G_1,\Sigma_1)$ be obtained from $(G,\Sigma)$ by deleting $v$,
and let $(H,\Omega)$ be obtained from $(G,\Sigma)$ by deleting $v$ and $w$. If $\Is(H,\Omega) = \I(H,\Omega)$, then
\begin{equation*}
\Is(G,\Sigma)  = \Is(G_1,\Sigma_1) + \mathbb{N}_1^2.
\end{equation*}
Furthermore, if $\I(G_1,\Sigma_1) = \Is(G_1,\Sigma_1)$, then $\I(G,\Sigma) = \Is(G,\Sigma)$.
\end{thm}
\begin{proof}
By Theorem~\ref{thm:contractingedge}, 
\begin{equation*}
\Is(G_1,\Sigma_1) + \mathbb{N}_1^2\subseteq \Is(G,\Sigma)\subseteq \I(G,\Sigma).
\end{equation*}
(As edge $e$, we can use the edge connecting $v$ and $w$.)

For the converse inclusion, suppose 
\begin{equation*}
(p,q)\in \Is(G,\Sigma)\setminus (\Is(G_1,\Sigma_1) + \mathbb{N}_1^2).
\end{equation*}
We may assume that $V(G) = \{1,\ldots,n\}$ and that $v=n$ and $w=n-1$. Let $R=V(G)\setminus\{v,w\}$. We claim that
any matrix $A=[a_{i,j}]\in S(G,\Sigma)$ with $p$ positive and $q$ negative eigenvalues, and the SAP has $a_{v,v} = 0$. To see this, suppose for a contradiction that $a_{v,v}\not=0$.  If $a_{v,v}>0$, then 
\begin{equation*}
B = A/A[v] =\begin{bmatrix}
A[R] & A[R,w]\\
A[w,R] & a_{w,w}-a_{w,v} a_{v,v}^{-1} a_{v,w}
\end{bmatrix}
\in S(G_1,\Sigma_1)
\end{equation*}
has $p-1$ positive and $q$ negative eigenvalues. Suppose $B$ does not have the SAP\@. Then there exists a nonzero symmetric matrix $X=[x_{i,j}]$ with $x_{i,j} = 0$ if $i$ and $j$ are adjacent in $(G_1,\Sigma)$ or if $i=j$, such that $B X = 0$. Let 
\begin{equation*}
Y = [y_{i,j}]=\begin{bmatrix}
X[R] & X[R,w] & -a_{v,v}^{-1} a_{v,w} X[R,w]\\
X[w,R] & 0 & 0\\
-a_{v,v}^{-1} a_{v,w} X[w,R] & 0 & 0
\end{bmatrix}.
\end{equation*}
Then $y_{i,j} = 0$ if $i$ and $j$ are adjacent in $(G,\Sigma)$ or if $i=j$, and $A Y = 0$. Hence $A$ would not satisfy the SAP\@. This contradiction shows that $B$ has the SAP\@. Hence $(p-1,q) \in \Is(G_1,\Sigma_1)$, and so $(p,q) \in \{(1,0)\}+\Is(G_1,\Sigma_1)$, which is a contradiction. If $a_{v,v} < 0$, then 
$B = A/A[v] \in S(G_1,\Sigma_1)$ has $p$ positive and $q-1$ negative eigenvalues. In the same way as above one can show that $B$ has the SAP\@. Then $(p,q) \in \{(0,1)\}+\Is(G_1,\Sigma_1)$, a contradiction. These contradictions show that $a_{v,v}=0$. 

A calculation shows that $A/A[\{v,w\}] = A[R]$.
Since $A[\{v,w\}]$ has one  positive and one negative eigenvalue, $A[R] = A/A[\{v,w\}]$ has $p-1$ positive and $q-1$ negative eigenvalues. Since $A[R]\in S(H,\Omega)$, $(p-1,q-1)\in \I(H,\Omega) = \Is(H,\Omega)$. By Theorem~\ref{thm:stableinertiaminor}, $\Is(H,\Omega)+\mathbb{N}_1^2\subseteq \Is(G_1,\Sigma_1)$, and so $(p,q-1)\in \Is(G_1,\Sigma_1)$. Hence $(p,q)\in \Is(G_1,\Sigma_1)+\{(0,1)\}\subseteq \Is(G_1,\Sigma_1)+\mathbb{N}_1^2$, which is a contradiction. 

Assume now furthermore that $\I(G_1,\Sigma_1) = \Is(G_1,\Sigma_1)$. Suppose $(p,q) \in \I(G,\Sigma)\setminus (\Is(G_1,\Sigma_1)+\mathbb{N}_1^2)$. Then any matrix $A=[a_{i,j}]\in S(G,\Sigma)$ with $p$ positive and $q$ negative eigenvalues has $a_{v,v}=0$. The proof of this is similar as above. Then $(p-1,q-1)\in \I(H,\Omega)=\Is(H,\Omega)$ and then in the same way as above $(p,q)\in \Is(G_1,\Sigma_1)+\mathbb{N}_1^2$.
\end{proof}

\begin{thm}\label{thm:supressingdegreetwo}
Let $(G,\Sigma)$ be a signed graph, let $v$ be a vertex of degree two in $G$ with two neighbors $u$ and $w$, and suppose the edges incident to $v$ are even. Let $(G_1,\Sigma)$ and $(G_2,\Sigma_2)$ be obtained from $(G,\Sigma)$ by deleting $v$ and adding between $u$ and $w$ respectively an even edge and an odd edge. Let $(H,\Omega)$ be obtained from $(G,\Sigma)$ by deleting $v$ and identifying $u$ and $w$.
If $\Is(H,\Omega) = \I(H,\Omega)$, then 
\begin{equation*}
\Is(G,\Sigma) = [\{(1,0)\}+\Is(G_1,\Sigma)]\cup [\{(0,1)\}+\Is(G_2,\Sigma_2)].
\end{equation*}
Furthermore, if $\Is(G_1,\Sigma) = \I(G_1,\Sigma)$ and $\Is(G_2,\Sigma_2) = \I(G_2,\Sigma_2)$, then $\I(G,\Sigma) = \Is(G,\Sigma)$
\end{thm}
\begin{proof}
By Theorem~\ref{thm:contractingedge}, 
\begin{equation*}
[\{(1,0)\}+\Is(G_1,\Sigma)]\cup [\{(0,1)\}+\Is(G_2,\Sigma_2)]\subseteq \Is(G,\Sigma).
\end{equation*}
(As edge $e$, we can use the edge connecting $v$ and $w$.)

For the converse inclusion, suppose $(p,q) \in \Is(G,\Sigma)$ and 
\begin{equation*}
(p,q)\not\in [\{(1,0)\}+\Is(G_1,\Sigma)]\cup [\{(0,1)\}+\Is(G_2,\Sigma_2)].
\end{equation*}
We may assume that $V(G) = \{1,\ldots,n\}$ and that $v=n$, $u=n-1$, and $w=n-2$.  Let $R=V(G)\setminus\{u,v,w\}$ and $S=\{u,w\}$. We claim that any matrix $A=[a_{i,j}]\in S(G,\Sigma)$ with $p$ positive and $q$ negative eigenvalues, and the SAP has $a_{v,v} = 0$. To see this, suppose for a contradiction that $a_{v,v}\not=0$. If $a_{v,v}>0$, then 
$B = A/A[v]  \in S(G_1,\Sigma)$
has $p-1$ positive and $q$ negative eigenvalues. 
Suppose $B$ does not have the SAP\@. Then there exists a nonzero symmetric matrix $X=[x_{i,j}]$ with $x_{i,j} = 0$ if $i$ and $j$ are adjacent in $(G_1,\Sigma)$ or if $i=j$, such that $B X = 0$. Let 
\begin{equation*}
Y = [y_{i,j}]=\begin{bmatrix}
X[R] & X[R,S] & -a_{v,v}^{-1} X[R,S]A[S,v]\\
X[S,R] & 0 & 0\\
-a_{v,v}^{-1} A[v,S] X[S,R] & 0 &0
\end{bmatrix}.
\end{equation*}
Then $y_{i,j} = 0$ if $i$ and $j$ are adjacent in $(G,\Sigma)$ or if $i=j$, and $A Y = 0$. Hence $A$ would not satisfy the SAP\@. This contradiction shows that $B$ has the SAP\@.
Hence $(p,q) \in [\{(1,0)\}+\Is(G_1,\Sigma)]$, a contradiction. If $a_{v,v} < 0$, then $B = A/A[v] \in S(G_2,\Sigma_2)$ has $p$ positive and $q-1$ negative eigenvalues, and in the same way as above one can show that $B$ has the SAP\@. Then $(p,q) \in [\{(0,1)\}+\Is(G_2,\Sigma_2)]$, a contradiction. These contradictions show that $a_{v,v}=0$.

We claim that $A/A[\{v,u\}]\in S(H,\Omega)$. By simultaneously scaling the $u$th row and column, we may assume that $a_{v,u} = -1$. Similarly, by simultaneously scaling the $w$th row and column we may assume that $a_{v,w} = -1$. Let $c = a_{w,w} + a_{u,u}  -2a_{u,w}$
Then 
\begin{equation*}
A/A[\{v,u\}] = \begin{bmatrix}
A[R] & A[R,u]-A[R,w]\\
A[u,R]-A[w,R] & c
\end{bmatrix}\in S(H,\Omega).
\end{equation*}
Since $A[\{v,u\}]$ has one  positive and one negative eigenvalue, $A/A[\{v,u\}]$ has $p-1$ positive and $q-1$ negative eigenvalues. Hence $(p-1,q-1)\in \I(H,\Omega) = \Is(H,\Omega)$. By Theorem~\ref{thm:stableinertiaminor}, $\{(1,0)\}+\Is(H,\Omega)\subseteq \Is(G_2,\Sigma_2)$
and $\{(0,1)\}+\Is(H,\Omega)\subseteq \Is(G_1,\Sigma)$, and so $(p,q-1)\in \Is(G_2,\Sigma_2)$ and $(p-1,q) \in \Is(G_1,\Sigma)$. Hence $(p,q)\in [\{(1,0)\}+\Is(G_1,\Sigma)]\cup [\{(0,1)\}+\Is(G_2,\Sigma_2)]$, which is a contradiction.

The proof that $\I(G,\Sigma) = \Is(G,\Sigma)$ if $\Is(G_1,\Sigma) = \I(G_1,\Sigma)$ and $\Is(G_2,\Sigma_2) = \I(G_2,\Sigma_2)$ goes along the same lines.
\end{proof}

By applying induction on $k$ and using Propositions~\ref{prop:Kn=}~and~\ref{prop:signedKn}, and Theorem~\ref{thm:supressingdegreetwo}, one can prove the following proposition.
\begin{prop}
If $k$ is an even positive integer, then 
\begin{equation*}
\Is(C_k^e) = \I(C_k^e) = \mathbb{N}^2_{[k-1,k]}\cup \{(2p+1,k-3-2p)~|~p=0,\ldots,k/2-2\}
\end{equation*} 
and 
\begin{equation*}
\Is(C_k^o) = \I(C_k^o) = \mathbb{N}^2_{[k-1,k]}\cup \{(2p,k-2-2p)~|~p=0.\ldots,k/2-1\}.
\end{equation*}
If $k$ is an odd integer $\geq 3$, then 
\begin{equation*}
\Is(C_k^e) = \I(C_k^e) = \mathbb{N}^2_{[k-1,k]}\cup \{(2p,k-2-2p)~|~p=0,\ldots,(k-3)/2\}
\end{equation*} 
and 
\begin{equation*}
\Is(C_k^o) = \I(C_k^o) = \mathbb{N}^2_{[k-1,k]}\cup \{(2p+1,k-3-2p)~|~p=0.\ldots,(k-3)/2\}.
\end{equation*}
\end{prop}
%

\begin{cor}
If $k$ is an integer $\geq 2$, then $\Is(C_k) = \mathbb{N}^2_{[k-2,k]}$.
\end{cor}

\begin{prop}\label{prop:signedK4oe}
$\Is(K_4^d) = \I(K_4^d) = \mathbb{N}_{[2,4]}$.
\end{prop}
\begin{proof}
We may assume that $K_4^d$ contains one odd edge $e$.
Since $K_4^d$ contains an odd triangle which can be obtained by deleting one vertex, $\mathbb{N}_{[2,4]}\subseteq \Is(K_4^d)$. Suppose for a contradiction that there exists a matrix $A\in S(K_4^d)$ with nullity $>2$. Let $f=vw$ be an even edge adjacent to $e$, where $v$ is incident to $e$. As $A$ has nullity $>2$, there exists a nonzero vector $x\in \ker(A)$ with $x_v = x_w = 0$. The other components of $x$ must be nonzero. Let $a_v$ and $a_w$ be the $v$th and $w$th row of $A$, respectively. From $a_w x = 0$ it follows that one component is negative and the other positive. From $a_v x = 0$ it follows that both are either negative or positive. This contradiction shows that $A$ has nullity at most two.
Hence $\Is(K_4^d) = \I(K_4^d) = \mathbb{N}_{[2,4]}$.
\end{proof}

\begin{prop}
$\Is(K_{2,3}^o)=\I(K_{2,3}^o) = \mathbb{N}_{[3,5]}$ and
\begin{equation*}
\Is(K_{2,3}^e)=\I(K_{2,3}^e) = [\{(1,1)\}+\mathbb{N}_{[0,3]}]\cup \{(0,4),(0,5),(4,0),(5,0) \}.
\end{equation*}
\end{prop}
\begin{proof}
Since $K_{2,3}^e$ is sign-equivalent to $(K_{2,3},E(K_{2,3}))$, $(p,q)\in \Is(K_{2,3}^e)$ if and only if $(q,p)\in \Is(K_{2,3}^e)$. A similar statement holds for $K_{2,3}^o$.

Since the signed graph $K_{2,3}^e$ is bipartite, $(3,0)\not\in \I(K_{2,3}^e)$ and $(4,0)\in \Is(K_{2,3}^e)$. Since $\mu(K_{2,3}) = 3$ (see \cite{CdeV2}), $(1,1) \in \Is(K_{2,3}^e)$. Hence $\I(K_{2,3}^e) = \Is(K_{2,3}^e)=(\{(1,1)\}+\mathbb{N}_{[0,3]})\cup \{(0,4),(0,5),(4,0),(5,0) \}$.

We now prove that $\I(K_{2,3}^o) = \Is(K_{2,3}^o)=\mathbb{N}_{[3,5]}$. By resigning if necessary, we may assume that $K_{2,3}^o$ has exactly one odd edge $e$. Removing the end of $e$ with degree two leaves an even $4$-cycle, and removing another vertex with degree two leaves an odd $4$-cycle. Since $\Is(C_4^o)\cup \Is(C_4^e) = \mathbb{N}_{[2,4]}$, we obtain from Theorem~\ref{thm:stableinertiaminor} that $\mathbb{N}_{[3,5]}\subseteq \Is(K_{2,3}^o)\subseteq \I(K_{2,3}^o)$.

Suppose that there exists a $(p,q)\in \I(K_{2,3}^o)$ with $p+q\leq 2$. Then there exists a matrix $A\in S(K_{2,3}^o)$ with nullity at least three. Let $v$ be the vertex of degree three which is incident to $e$. Let $w$ be the other vertex of degree three, and let $u_1,u_2,u_3$ be the vertices of degree two, where we assume that $u_1$ is incident to $e$. Since $A$ has nullity at least three, there exists a nonzero vector $x\in \ker(A)$ with $x_v = x_{u_2} = 0$.  For each vertex $i$ of $K_{2,3}^o$, let $a_i$ denote the $i$th row of $A$. From $a_{u_2}x = 0$, we obtain that $x_w = 0$. From $a_v x = 0$, we obtain that $x_{u_1}$ and $x_{u_3}$ must have the same sign. However, from $a_w x = 0$, we obtain that $x_{u_1}$ and $x_{u_3}$ must have different signs. This contradiction shows that $\Is(K_{2,3}^o) = \I(K_{2,3}^o)=\mathbb{N}_{[3,5]}$.
\end{proof}

\begin{cor}
$\Is(K_{2,3}) = \mathbb{N}^2_{[3,5]}\cup \{(1,1)\}$.
\end{cor}

\begin{prop}
$\Is(K_{3,3}) = \I(K_{3,3}) = \{(1,1)\}\cup \oN^2_{[3,6]}$.
\end{prop}
\begin{proof}
Since $K_{3,3}$ has a $K_4$-minor, $\oN^2_{[3,6]}\subseteq \I_s(K_{3,3})$ by Corollary~\ref{cor:graphminormonotonicity} and Corollary~\ref{cor:IsKn}. Since $\xi(K_{3,3}) = 4$ and $\mr_+(K_{3,3})\geq 3$ (see \cite{MR2111528}), $(1,1)\in \Is(K_{3,3})$. Hence $\{(1,1)\}\cup \oN^2_{[3,6]} \subseteq \Is(K_{3,3})$. Since $\mr(K_{3,3})=2$, $\mr_+(K_{3,3})=3$, and $\Is(K_{3,3})\subseteq \I(K_{3,3})$, $\Is(K_{3,3}) = \I(K_{3,3}) = \{(1,1)\}\cup \oN^2_{[3,6]}$.
\end{proof}

\section{Disjoint unions and $1$-sums}\label{sec:01splits}

A \emph{cut-vertex} of a graph is a vertex whose deletion
increases the number of connected components.  
A connected graph $G$ is \emph{{\rm 2}-connected} if $G$ has at least three vertices and no cut-vertex.  A \emph{block} of a graph is a maximal connected
subgraph without a cut-vertex. 
Let $G$ be a graph and let $C$ be a block of $G$.
The \emph{thin out} of $C$ in $G$ is the graph obtained from $C$ by adding a pendant edge
to each cut vertex $v$ of $G$ contained in $C$. So the thin out of 
$G$ is a subgraph of $G$. These definitions extend to signed graphs. In this section we show how the stable inertia set of a signed graph $(G,\Sigma)$ can be determined from the stable inertia sets of the thin out of the $2$-connected blocks in $(G,\Sigma)$. A similar statement holds for graphs.

The following lemma will be used in the proofs of the theorems of this section.
\begin{lemma}\label{lem:SAPnonsing}
Let $(G,\Sigma)$ be a signed graph and let $A\in S(G,\Sigma_1)$. Suppose $S_1$ and $S_2$ are nonempty disjoint subsets of $V(G)$ such that $A[S_1,S_2] = 0$. Let $R = V(G)\setminus (S_1\cup S_2)$. If $A$ has the SAP, then it is not possible that there are nonzero vectors $y$ and $z$ such that $A[S_1\cup R,S_1]y = 0$ and $A[S_2\cup R,S_2]z = 0$.
\end{lemma}
\begin{proof}
We may assume that $S_1=\{1,\ldots,k\}$, $R=\{k+1,\ldots, m\}$ and $S=\{m+1,\ldots, n\}$. If $A[S_1\cup R,S_1]y = 0$ and $A[S_2\cup R,S_2]z = 0$, then the matrix
\begin{equation*}
X = [x_{i,j}] = \begin{bmatrix}
0 & 0 & y z^T\\
0 & 0 & 0\\
z y^T & 0 & 0
\end{bmatrix}
\end{equation*}
is nonzero, satisfies $x_{i,j} = 0$ if $i$ and $j$ are adjacent in $G$ or if $i=j$, and $A X = 0$, which is a contradiction.
\end{proof}

The \emph{disjoint union} of two signed graphs $(G_1,\Sigma_1)$ and $(G_2,\Sigma_2)$ with disjoint vertex sets is the signed graph $(G_1\cup G_2,\Sigma_1\cup \Sigma_2)$.

\begin{thm}\label{thm:Is0split}
Suppose $(G,\Sigma)$ is the disjoint union of $(G_1,\Sigma_1)$ and $(G_2,\Sigma_2)$, and let $n_1 = \lvert V_1\rvert$ and $n_2=\lvert V_2\rvert$. Then $\Is(G,\Sigma) = [\Is(G_1,\Sigma_1)+\mathbb{N}^2_{n_2}]\cup [\Is(G_2,\Sigma_2)+\mathbb{N}^2_{n_1}]$.
\end{thm}
\begin{proof}
Since $(G_1,\Sigma_1)$ and $(G_2,\Sigma_2)$ are subgraphs of $(G,\Sigma)$, $\Is(G_1,\Sigma_1)+ \mathbb{N}^2_{n_2} \subseteq \Is(G,\Sigma)$ and $\Is(G_2,\Sigma_2)+\mathbb{N}^2_{n_1}\subseteq \Is(G,\Sigma)$.
Hence $[\Is(G_1,\Sigma_1)+\mathbb{N}^2_{n_2}]\cup [\Is(G_2,\Sigma_2)+\mathbb{N}^2_{n_1}]\subseteq \Is(G,\Sigma)$.

To see the converse inclusion, let $(p,q)\in \Is(G,\Sigma)$. Let $A\in \mathcal{S}(G,\Sigma)$ have the SAP, $p$ positive and $q$ negative eigenvalues. Since $A$ has the SAP, $A[V(G_1)]$ or $A[V(G_2)]$ is nonsingular by Lemma~\ref{lem:SAPnonsing}. We may assume that $A[V(G_1)]$ is nonsingular, as the other case is similar.

Suppose that $A[V(G_2)]$ does not have the SAP\@. Then there exists a nonzero symmetric matrix $X=[x_{i,j}]$ with $x_{i,i}=0$ for all $i\in V(G_2)$, $x_{i,j}=0$ for all $ij \in E(G_2)$, and $A[V(G_2)] X = 0$. Let
\begin{equation*}
Y = [y_{i,j}] = \begin{bmatrix}
0 & 0\\
0 & X
\end{bmatrix}.
\end{equation*}
Then $y_{i,i}=0$ for all $i\in V(G)$, $y_{i,j}= 0$ for all $ij\in E(G)$, and $A Y = 0$. Hence $A$ does not have the SAP, contradicting the assumption.

Thus, $(p,q)\in \Is(G_2,\Sigma_2) + \mathbb{N}^2_{n_1}$.
\end{proof}

From Theorem~\ref{thm:Is0split}, we immediately obtain the following corollaries.

\begin{cor}\label{cor:graphdisjointunion}
Let $G$ be a disjoint union of $G_1$ and $G_2$. Let $n_1 = |G_1|$ and $n_2=|G_2|$.
Then $\Is(G) = [\Is(G_1)+\mathbb{N}^2_{n_2}]\cup [\Is(G_2)+\mathbb{N}^2_{n_1}]$.

\end{cor}

\begin{cor}\label{cor:nu0split}
If $(G,\Sigma)$ is the disjoint union of $(G_1,\Sigma_1)$ and $(G_2,\Sigma_2)$, then $\nu(G,\Sigma) = \max\{\nu(G_1,\Sigma_1), \nu(G_2,\Sigma_2)\}$.
\end{cor}

\begin{cor}\label{cor:xi0split}
If $(G,\Sigma)$ is the disjoint union of $(G_1,\Sigma_1)$ and $(G_2,\Sigma_2)$, then $\xi(G,\Sigma) = \max\{\xi(G_1,\Sigma_1), \xi(G_2,\Sigma_2)\}$.
\end{cor}

If $G_1 = (V_1,E_1)$ and $G_2=(V_2,E_2)$ are subgraphs of $G$ such that $G=G_1\cup G_2$ and $V_1\cap V_2=\{v\}$, then $G$ is called the \emph{$1$-sum} of $G_1$ and $G_2$ at $v$. 
Let $(G,\Sigma)$ be a signed graph. If $G$ is the $1$-sum of $G_1$ and $G_2$ at $v$, then $(G,\Sigma)$ is called the \emph{$1$-sum} of $(G_1,\Sigma\cap E_1)$ and $(G_2,\Sigma\cap E_2)$ at $v$. 

For $(p_1,q_1), (p_2,q_2) \in \mathbb{N}^2$, define $(p_1,q_1) + (p_2,q_2) = (p_1+p_2,q_1+q_2)$.

\begin{thm}\label{thm:Is1split}
Let $(G,\Sigma)$ be a connected signed graph and suppose $(G,\Sigma)$ is the $1$-sum of $(G_1,\Sigma_1)$ and $(G_2,\Sigma_2)$ at $v$, with both $E(G_1)$ and $E(G_2)$ nonempty. Let $n_1$ and $n_2$ be the number of vertices of $G_1$ and $G_2$, respectively. For $i=1,2$, let $(H_i,\Sigma_i)$ be the signed graph obtained from $K_2$ and $G_i$ by identifying a vertex of $K_2$ with $v$. 
Then $\Is(G,\Sigma) = [\Is(H_1,\Sigma_1)+\oN^2_{n_2-2}]\cup [\Is(H_2,\Sigma_2)+\oN^2_{n_1-2}]$.
\end{thm}
\begin{proof}
By negating around vertices, we may assume that for $i=1,2$ at least one edge of $G_i$ incident to $v$ is even.
Since $(H_1,\Sigma_1)$ and $(H_2,\Sigma_2)$ are subgraphs of $(G,\Sigma)$, $\Is(H_1,\Sigma_1)+\oN_{n_2-2}\subseteq \Is(G,\Sigma)$ and $\Is(H_2,\Sigma_2)+\oN_{n_1-2}\subseteq \Is(G,\Sigma)$. Hence $[\Is(H_1,\Sigma_1)+\oN_{n_2-2}]\cup [\Is(H_2,\Sigma_2)+\oN_{n_1-2}]\subseteq \Is(G,\Sigma)$.

To see the converse inclusion, let $A\in \mathcal{S}(G,\Sigma)$ have the SAP\@. For $i=1,2$, let $V_i = V(G_i)$, $E_i = E(G_i)$, and $S_i=V_i\setminus\{v\}$.

Suppose $A[S_1]$ is nonsingular. We may write
\begin{equation*}
A  = \begin{bmatrix}
     A[S_1] & A[S_1,v] & 0  \\
     A[v,S_1] & a_{v,v} & A[v,S_2] \\
     0 & A[S_2,v] & A[S_2]
   \end{bmatrix}.
\end{equation*}
Let
\begin{equation*}
P = \begin{bmatrix}
   I & -A[S_1]^{-1} A[S_1,v] & 0  \\
   0 & 1 & 0  \\
   0 & 0 & I  \\
\end{bmatrix}.
\end{equation*}
By Sylvester's Law of Inertia,
\begin{equation*}
P^T AP = \begin{bmatrix}
   A[S_1] & 0 & 0  \\
   0 & a_{v,v}-A[v,S_1]A[S_1]^{-1} A[S_1 ,v] & A[v,S_2 ]  \\
   0 & A[S_2,v] & A[S_2]
\end{bmatrix}
\end{equation*}
has the same inertia as $A$. The matrix
\begin{equation*}
B := \begin{bmatrix}
a_{v,v}-A[v,S_1]A[S_1]^{-1}A[S_1 ,v] & A[v,S_2 ]  \\
A[S_2 ,v] & A[S_2 ]
    \end{bmatrix}
\end{equation*}
belongs to $\mathcal{S}(G_2,\Sigma_2)$ and $\pin(B) + \pin(A[S_1]) = \pin(A)$.

Suppose for a contradiction that $B$ does not have the SAP\@. Then there exists a nonzero symmetric matrix $X=[x_{i,j}]$ with $x_{i,i}=0$ for all $i\in V_2$, $x_{i,j}=0$ for all $ij \in E_2$, and $B X = 0$.
Let $Z = -A[S_1]^{-1} A[S_1,v]X[v,S_2]$ and
\begin{equation*}
Y = [y_{i,j}] = \begin{bmatrix}
0 & 0 & Z\\
0 & 0 & X[v,S_2]\\
Z^T & X[S_2,v] & X[S_2]
\end{bmatrix}.
\end{equation*}
Then $y_{i,i}=0$ for all $i\in V$, $y_{i,j}= 0$ for all $ij\in E$, and $A Y = 0$. Hence $A$ does not have the SAP, contradicting the assumption.

Hence $\pin(A) \in \Is(G_2,\Sigma_2)+\oN^2_{n_1-1}$ and since $G_2$ is a subgraph of $H_2$, $\pin(A)\in \Is(H_2,\Sigma_2)+\oN_{n_1-2}$. The case where $A[S_2]$ is nonsingular can be done similarly.

Therefore we may assume that both $A[S_1]$ and $A[S_2]$ are singular.
Since $A$ has the SAP, it is, by Lemma~\ref{lem:SAPnonsing}, not possible that there are nonzero vectors $y$ and $z$ such that $A[V_1,S_1]y=0$ and $A[V_2,S_2]z=0$; say there is no nonzero vector $y$ such that $A[V_1,S_1]y=0$.

Since there is no nonzero $y$ with $A[S_1\cup\{v\},S_1]y=0$, $A[S_1]$ has nullity $1$. Let $x\in \ker(A[S_1])$ be nonzero. Let $w\in S_1$ with $x_w\not=0$. Let $Q=S_1\setminus\{w\}$. Then $A[Q]$ is nonsingular, because if $A[Q]$ were singular another (independent) vector could be constructed in $\ker(A[S_1])$. We may write
\begin{equation*}
A = \begin{bmatrix}
A[Q] & A[Q,w] & A[Q,v] & 0\\
A[w,Q] & a_{w,w} & a_{w,v} & 0\\
A[v,Q] & a_{v,w} & a_{v,v} & A[v,S_2]\\
0 & 0 & A[S_2,v] & A[S_2]
\end{bmatrix}.
\end{equation*}
Let
\begin{equation*}
P=\begin{bmatrix}
I & -A[Q]^{-1} A[Q,\{w,v\}] & 0\\
0 & I_2 & 0\\
0 & 0 & I
\end{bmatrix}.
\end{equation*}
Here $I_2$ denotes the $2\times 2$ identity matrix.
By Sylvester's law of Inertia
\begin{equation*}
P^T A P = \begin{bmatrix}
A[Q] & 0 & 0 & 0\\
0 & b_{w,w} & b_{w,v} & 0\\
0 & b_{v,w} & b_{v,v} & A[v,S_2]\\
0 & 0 & A[S_2,v] & A[S_2]
\end{bmatrix},
\end{equation*}
where
\begin{equation*}
\begin{bmatrix}
b_{w,w} & b_{w,v}\\
b_{v,w} & b_{v,v}
\end{bmatrix} = A[\{v,w\}] - A[\{v,w\},Q] A[Q]^{-1} A[Q,\{v,w\}],
\end{equation*}
has the same inertia as $A$.
Let
\begin{equation*}
B := \begin{bmatrix}
A[\{v,w\}] - A[\{v,w\},Q] A[Q]^{-1} A[Q,\{v,w\}] & A[\{v,w\},S_2]\\
A[S_2,\{v,w\}] & A[S_2]
\end{bmatrix}.
\end{equation*}
The matrix $B$ satisfies $\pin(B) + \pin(A[Q]) = \pin(A)$.
Since $A[S_1]$ is singular, we know that $b_{w,w} = a_{w,w}-A[w,Q] A[Q]^{-1} A[Q,w] = 0$.
Suppose $b_{v,w} = a_{v,w} -A[v,Q] A[Q]^{-1} A[Q,w] = 0$. Then the vector
\begin{equation*}
a = \begin{bmatrix}
-A[Q]^{-1} A[Q,w]\\
1
\end{bmatrix}
\end{equation*}
belongs to $\ker(A[V_1,S_1])$, contradicting the assumption. Therefore $b_{v,w}\not=0$, that is, $B\in\mathcal{S}(H_2,\Sigma_2)$.

Suppose for a contradiction that $B$ does not have the SAP\@. Then there is a nonzero symmetric matrix $X=[x_{i,j}]$ such that $x_{i,j}=0$ if $i=j$ or $ij\in E$, and $B X = 0$.
Let
\begin{equation*}
Z = -A[Q]^{-1} A[Q,\{v,w\}]X[\{v,w\},S]
\end{equation*}
and
\begin{equation*}
Y = [y_{i,j}] = \begin{bmatrix}
0 & 0 & Z\\
0 & 0 & X[\{v,w\},S_]\\
Z^T & X[S_1,\{v,w\}] & X[S_1]
\end{bmatrix}.
\end{equation*}
Then $Y$ is nonzero, $y_{i,j}=0$ if $i=j$ or $ij\in E$, and $A Y = 0$. Hence $A$ would not have the SAP if $B$ did not.
Hence $\pin(A)\in \Is(H_2,\Sigma_2)+\oN_{n_1-2}$.
\end{proof}

\begin{cor}\label{cor:Is1sumgraph}
Let $G$ be a connected graph and suppose $G$ is the $1$-sum of $G_1$ and $G_2$ at $v$, both containing at least one edge. Let $n_1$ and $n_2$ be the number of vertices of $G_1$ and $G_2$, respectively. For $i=1,2$, let $H_i$ be the graph obtained from $K_2$ and $G_i$ by identifying a vertex of $K_2$ with $v$. 
Then $\Is(G) = [\Is(H_1)+\oN^2_{n_2-2}]\cup [\Is(H_2)+\oN^2_{n_1-2}]$.
\end{cor}

We now show how the stable inertia set of a signed graph $(G,\Sigma)$ can be determined from stable inertia set of the thin out of each $2$-connected block in $(G,\Sigma)$. We do this by induction on the number of vertices in $G$. Suppose $(G,\Sigma)$ is not the thin out of a $2$-connected block. If $(G,\Sigma)$
is a disjoint union of $(G_1,\Sigma_1)$ and $(G_2,\Sigma_2)$, then we apply Theorem~\ref{thm:Is0split} on $(G_1,\Sigma_1)$ and $(G_2,\Sigma)$. We may therefore assume that $(G,\Sigma)$ is connected. If $(G,\Sigma)$ has a $2$-connected block, then 
we apply Theorem~\ref{thm:Is1split}.  If $(G,\Sigma)$ has no $2$-connected block, then each block is either a single edge or a class of parallel edges. We may assume that if a block has parallel edges, then the block is equal to $K_2^=$.
If $(G,\Sigma)$ is a $1$-sum of $(G_1,\Sigma\cap E(G_1))$ and $(G_2,\Sigma\cap E(G_2))$ with both $G_1$ and $G_2$ containing at least two edges, then we apply Theorem~\ref{thm:Is1split}. Hence, we may assume that if $(G,\Sigma)$ is a $1$-sum of $(G_1,\Sigma\cap E(G_1))$ and $(G_2,\Sigma\cap E(G_2))$, then $G_1$ or $G_2$ contains at most one edge.

Suppose now that $(G,\Sigma)$ has a block that is equal to $K_2^=$. Then the number, $n$, of vertices  of $G$ is at most four and $\Is(K_2^=) + \mathbb{N}^2_{n-2} = \mathbb{N}^2_{[n-2,n]}\subseteq \Is(G,\Sigma)$. If $(p,q) \in \Is(G,\Sigma)$ and $(p,q)\not\in \mathbb{N}^2_{[n-2,n]}$, then $\xi(G,\Sigma) \geq 3$. Hence there exists a matrix $A\in S(G,\Sigma)$ with $\nullity A \geq 3$. Let $u,v$ be the vertices of $K_2^=$. Since $\nullity A\geq 3$, there exists a nonzero vector $x$ with $x_u = x_v = 0$. However, this forces $x_w=0$ for any other vertex $w$ of $(G,\Sigma)$. This contradiction shows that $\mathbb{N}^2_{[n-2,n]} = \Is(G,\Sigma)$.

Suppose finally that $(G,\Sigma)$ has no blocks that are equal to $K_2^=$. Then $(G,\Sigma)$ is sign-equivalent to $(T,\emptyset)$, where $T$ is a tree. We then use Proposition~\ref{prop:treestabinertia}.

\section{Some characterizations}\label{sec:lowvalues}

A signed graph has no $K_2^=$-minor if and only if it has no odd cycles, that is, if it is bipartite. Since $\nu(K_2^=)=2$, each signed graph $(G,\Sigma)$ with $\nu(G,\Sigma)\leq 1$ has no $K_2^=$-minor, that is, is bipartite. In fact, $K_2^=$ is the only obstruction to having $\nu(G,\Sigma)\leq 1$.

\begin{thm}\label{thm:snu1}
A signed graph $(G,\Sigma)$ has $\nu(G,\Sigma)\leq 1$ if and only if $(G,\Sigma)$ is bipartite.
\end{thm}
\begin{proof}
If $\nu(G,\Sigma)\leq 1$, then $(G,\Sigma)$ has no $K_2^=$-minor, that is, $(G,\Sigma)$ has no odd cycle. Hence $(G,\Sigma)$ is bipartite.

For the converse, suppose $(G,\Sigma)$ is bipartite and $\nu(G,\Sigma)\geq 2$. Then, by Corollary~\ref{cor:nu0split}, there exists a component $(H,\Sigma_1)$ of $(G,\Sigma)$ with $\nu(H,\Sigma_1)\geq 2$. Hence there exists a positive semidefinite matrix $A\in S(H,\Sigma_1)$ with $\nullity(A)\geq 2$.
This contradicts Lemma~\ref{lem:smnull1}.
\end{proof}

From Theorem~\ref{thm:snu1} we obtain the following corollary.

\begin{cor}\cite{CdeV3}
A graph $G$ has $\nu(G)\leq 1$ if and only if $G$ is a forest.
\end{cor}
\begin{proof}
If $G$ is a forest, then for every subset $\Sigma\subseteq E(G)$, $(G,\Sigma)$ is bipartite. Hence $\nu(G,\Sigma)\leq 1$ for every subset $\Sigma\subseteq E(G)$. Since $\nu(G)=\max\{\nu(G,\Sigma)\mid \Sigma\subseteq E(G)\}$, we obtain that $\nu(G)\leq 1$.

For the converse, suppose $\nu(G)\leq 1$ and $G$ has a cycle $C$. Let $e$ be an edge of $C$. Then $\nu(G,\{e\})\leq \nu(G)\leq 1$. Since $(G,\{e\})$ has an odd cycle, we obtain a contradiction. 
\end{proof}

\begin{thm}\label{thm:M+1}
A signed graph $(G,\Sigma)$ has $M_+(G,\Sigma)\leq 1$ if and only if $(G,\Sigma)$ is connected bipartite.
\end{thm}
\begin{proof}
If $(G,\Sigma)$ has $M_+(G,\Sigma)\leq 1$, then $(G,\Sigma)$ is connected. Furthermore, since $\nu(G,\Sigma)\leq M_+(G,\Sigma)\leq 1$, $(G,\Sigma)$ is bipartite.

Conversely, if $(G,\Sigma)$ is connected bipartite, then, by Lemma~\ref{lem:smnull1}, each matrix $A\in S(G,\Sigma)$ has $\nullity(A)\leq 1$.
\end{proof}

From Theorem~\ref{thm:M+1}, we obtain the following corollary.

\begin{cor}
A graph $G$ has $M_+(G)\leq 1$ if and only if $G$ is a tree.
\end{cor}

We now characterize the signed graphs $(G,\Sigma)$ for which $\xi(G,\Sigma)\leq 1$.

\begin{thm}\label{thm:sxi1}
A signed graph $(G,\Sigma)$ has $\xi(G,\Sigma)\leq 1$ if and only if $(G,\Sigma)$ is sign-equivalent to the signed graph $(H,\emptyset)$, where $H$ is a graph whose underlying simple graph is a disjoint union of paths.
\end{thm}
\begin{proof}
If $\xi(G,\Sigma)\leq 1$, then $(G,\Sigma)$ has no $K_2^=$-minor, that is, $(G,\Sigma)$ has no odd cycle. Hence $(G,\Sigma)$ is bipartite. Furthermore, $(G,\Sigma)$ has no $K_{1,3}$-minor, that is, $G$ has no vertices with more than three neighbors. Thus $G$ is up to parallel edges a disjoint union of paths.

For the converse, suppose $G$ is a disjoint union of paths and $\xi(G,\Sigma)\geq 2$. Then, by Corollary~\ref{cor:xi0split}, there exists a component $(H,\Sigma_1)$ of $(G,\Sigma)$ with $\xi(H,\Sigma_1)\geq 2$. Since $M(H)\geq \xi(H,\Sigma_1)$, we obtain that $M(H)\geq 2$. Since $M(P)\leq 1$ for any path $P$, we obtain a contradiction.
\end{proof}

\begin{thm}
A signed graph $(G,\Sigma)$ has $M(G,\Sigma)\leq 1$ if and only if 
$(G,\Sigma)$ is sign-equivalent to a signed graph $(H,\emptyset)$, where $H$ is a graph whose underlying simple graph is a path.
\end{thm}
\begin{proof}
If $(G,\Sigma)$ has $M(G,\Sigma)\leq 1$, then $(G,\Sigma)$ is connected. Furthermore, $(G,\Sigma)$ has no $K_2^=$-minor and no $K_{1,3}$-minor, as $M(G,\Sigma)\geq \xi(G,\Sigma)\geq \xi(K_2^=) = 2$ and $M(G,\Sigma)\geq \xi(K_{1,3}) = 2$. Hence  $(G,\Sigma)$ is sign-equivalent to the signed graph $(H,\emptyset)$, where $H$ is a graph whose underlying simple graph is a path.

For the converse, suppose $(G,\Sigma)$ is signed-equivalent to a signed graph $(H,\emptyset)$, where $H$ is a graph whose underlying simple graph is a path. Since $M(H)\leq 1$, we obtain that $M(H,\emptyset)\leq 1$. Hence $M(G,\Sigma)\leq 1$.
\end{proof}

\bibliographystyle{plain}
\bibliography{../../../biblio}

\begin{thebibliography}{10}

\bibitem{BarFalHog2005a}
F.~Barioli, S.~Fallat, and L.~Hogben.
\newblock A variant on the graph parameters of {C}olin de {V}erdi\`ere:
  implications to the minimum rank of graphs.
\newblock {\em Electron. J. Linear Algebra}, 13:387--404 (electronic), 2005.

\bibitem{MR2547901}
W.~Barrett, H.~T. Hall, and R.~Loewy.
\newblock The inverse inertia problem for graphs: cut vertices, trees, and a
  counterexample.
\newblock {\em Linear Algebra Appl.}, 431(8):1147--1191, 2009.

\bibitem{MR2111528}
W.~Barrett, H.~van~der Holst, and R.~Loewy.
\newblock Graphs whose minimal rank is two.
\newblock {\em Electron. J. Linear Algebra}, 11:258--280 (electronic), 2004.

\bibitem{CdeV1}
Y.~Colin~de Verdi\`ere.
\newblock Sur un nouvel invariant des graphes et un crit\`ere de planarit\'e.
\newblock {\em J. Comb. Theory, Ser. B.}, 50:11--21, 1990.

\bibitem{CdeV2}
Y.~Colin~de Verdi\`ere.
\newblock On a new graph invariant and a criterion of planarity.
\newblock In N.~Robertson and P.~Seymour, editors, {\em Graph Structure
  Theory}, volume 147 of {\em Contemporary Mathematics}, pages 137--147.
  American Mathematical Society, Providence, Rhode Island, 1993.

\bibitem{CdeV3}
Y.~Colin~de Verdi\`ere.
\newblock Multiplicities of eigenvalues and tree-width of graphs.
\newblock {\em J. Comb. Theory, Ser. B.}, 74(2):121--146, 1998.

\bibitem{Diestel}
R.~Diestel.
\newblock {\em Graph Theory}.
\newblock Springer-Verlag, New York, second edition, 2000.

\bibitem{FalHog2007}
S.M. Fallat and L.~Hogben.
\newblock The minimum rank of symmetric matrices described by a graph: A
  survey.
\newblock {\em Linear Algebra Appl.}, 426(2--3):558--582, 2007.

\bibitem{MR0244285}
M.~Fiedler.
\newblock A characterization of tridiagonal matrices.
\newblock {\em Linear Algebra and Appl.}, 2:191--197, 1969.

\bibitem{gast.26511419740101}
M.~W. Hirsch and S.~Smale.
\newblock {\em Differential equations, dynamical systems, and linear algebra}.
\newblock Pure and applied mathematics; a series of monographs and textbooks,
  v. 60. New York, Academic Press [1974], 1974.

\bibitem{HvdH2007a}
L.~Hogben and H.~van~der Holst.
\newblock Forbidden minors for the class of graphs with $\xi({G}) \leq 2$.
\newblock {\em Linear Algebra Appl.}, 423:42--52, 2007.

\bibitem{Holst96a}
H.~van~der Holst.
\newblock {\em Topological and {S}pectral {G}raph {C}haracterizations}.
\newblock PhD thesis, University of Amsterdam, 1996.

\bibitem{Holst97a}
H.~van~der Holst.
\newblock Graphs with magnetic {S}chr{\"o}dinger operators of low corank.
\newblock {\em J. Comb. Theory, Ser. B.}, 84:311--339, 2002.

\bibitem{Holst2007b}
H.~van~der Holst.
\newblock The maximum corank of graphs with a $2$-separation.
\newblock {\em Linear Algebra Appl.}, 428(7):1587--1600, 2008.

\bibitem{HLovaszS2}
H.~van~der Holst, L.~Lov\'asz, and A.~Schrijver.
\newblock The {C}olin de {V}erdi\`ere graph parameter.
\newblock In {\em Graph Theory and Combinatorial Biology}, number~7 in
  Mathematical Studies, pages 29--85. Bolyai Society, 1999.

\bibitem{Kotlov2000a}
A.~Kotlov.
\newblock Spectral characterization of tree-width-two graphs.
\newblock {\em Combinatorica}, 20(1):147--152, 2000.

\bibitem{LovaszSchrijver98a}
L.~Lov\'{a}sz and A.~Schrijver.
\newblock A {B}orsuk theorem for antipodal links and a spectral
  characterization of linklessly embeddable graphs.
\newblock {\em Proceedings of the American Mathematical Society},
  126(5):1275--1285, 1998.

\bibitem{MR676405}
T.~Zaslavsky.
\newblock Signed graphs.
\newblock {\em Discrete Appl. Math.}, 4(1):47--74, 1982.

\end{thebibliography}

\end{document}